# Correspondence between Row-Column Determinants and Quasideterminants of Matrices over Quaternion Algebra

Aleks Kleyn and Ivan Kyrchei

Abstract. In this paper, we considered the theory of quasideterminants and row and column determinants. We considered the application of this theory to the solving of a system of linear equations in quaternion algebra. We established correspondence between row and column determinants and quasideterminants of matrix over quaternion algebra.

Contents

## 1. Preface

Linear algebra is a powerful tool that we use in different areas of mathematics, including calculus, analytic and differential geometry, theory of differential equations, optimal control theory. Linear algebra has accumulated a rich set of different methods. Since some methods have a common final result, this gives us the opportunity to choose the most effective method, depending on the nature of calculations.

At transition from linear algebra over a field to linear algebra over division ring, we want to save as much as possible tools that we regularly use. Already in the

Aleks_Kleyn@MailAPS.org.

kyrchei@lms.lviv.ua.

Pidstrygach Institute for Applied Problems of Mechanics and Mathematics, str.Naukova 3b, Lviv, Ukraine, 79005.

early XX century, shortly after Hamilton created quaternion algebra, mathematicians began to search the answer how looks like the algebra with noncommutative multiplication. In particular, there is a problem how to determine a determinant of a matrix with elements belonging to noncommutative ring. Such determinant is also called noncommutative determinant.

There were a lot of approaches to the definition of the noncommutative determinant. However none of the introduced noncommutative determinants maintained all those properties that determinant possessed for matrices over a field. Moreover, in paper [11], J. Fan proved that there is no unique definition of determinant which would expands the definition of determinant of real matrices for matrices over division ring of quaternions. Therefore, search for a solution of the problem to define a noncommutative determinant is still going on.

In this paper, we consider two approaches to define noncommutative determinant. Namely, we explore row-column determinants and quasideterminant.

Row-column determinants are an extension of the classical definition of the determinant, however we assume predetermined order of elements in each of the terms of the determinant. Using row-column determinants, we obtain the solution of system of linear equations over quaternion algebra according to Cramer's rule.

Quasideterminant appeared from the analysis of the procedure of matrix inversion. Using quasideterminant, the procedure of solving of a system of linear equations over quaternion algebra is similar to the method of Gauss.

There is common in definition of row and column determinants and quasideterminant. In both cases, we have not one determinant in correspondence to quadratic matrix of order $n$ with noncommutative entries, but certain set (there are $n^2$ quasideterminant, $n$ row determinants, and $n$ column determinants).

Today there is wide application of quasideterminants in linear algebra ([18, 21]), and in physics ([12, 13, 20]). Row and column determinants ([14]) introduced relatively recently are less well known. Purpose of this paper is establishment of the correspondence between row and column determinants and quasideterminants of matrix over quaternion algebra. The authors are hopeful that the establishment of this correspondence can provide mutual development of both theory quasideterminants and theory of row and column determinants.

## 2. Convention about Notations

There are different forms to write elements of matrix. In this paper, we denote $a_{ij}$ an element of the matrix $\mathbf{A}$. Index $i$ labels rows, and index $j$ labels columns.

We use the following notation for different minors of the matrix $\mathbf{A}$

- $\mathbf{a}_{i\,.}$ **:** the row with index $i$
- $\mathbf{A}_{S\,.}$ **:** the minor obtained from $A$ by selecting rows with index from the set $S$
- $\mathbf{A}^{i\,.}$ **:** the minor obtained from $A$ by deleting row $\mathbf{a}_{i\,.}$
- $\mathbf{A}^{S\,.}$ **:** the minor obtained from $A$ by deleting rows with index from the set $S$
- $\mathbf{a}_{.\,j}$ **:** column with index $j$
- $\mathbf{A}_{.\,T}$ **:** the minor obtained from $A$ by selecting columns with index from the set $T$
- $\mathbf{A}^{.\,j}$ **:** the minor obtained from $A$ by deleting column $\mathbf{a}_{.\,j}$

$\mathbf{A}^{\cdot T}$ **:** the minor obtained from $A$ by deleting columns with index from the set $T$

$\mathbf{A}_{.j}(\mathbf{b})$ **:** the matrix obtained from matrix $\mathbf{A}$ by replacing its column with number $j$ by column $\mathbf{b}$

$\mathbf{A}_{i.}(\mathbf{b})$ **:** the matrix obtained from matrix $\mathbf{A}$ by replacing its row with number $i$ by row $\mathbf{b}$

Considered notations can be combined. For instance, the record

$$\mathbf{A}^{ii}_{k.}(\mathbf{b})$$

means replacing of the row with number $k$ by matrix $\mathbf{b}$ followed by removal of the row with number $i$ and by removal of the column with number $i$.

As was noted in section [17]-2.2, we can define two types of matrix products: either product of rows of first matrix over columns of second one, or product of columns of first matrix over rows of second one. However according to the theorem [17]-2.2.5 this product is symmetric relative operation of transposition. Hence in the paper, we will restrict ourselves by traditional product of rows of first matrix over columns of second one; and we do not indicate clearly the operation like it was done in [17].

## 3. Preliminaries

Theory of determinants of matrices with noncommutative elements can be divided into three groups regarding their methods of definition. Denote $\mathrm{M}(n,\mathbf{K})$ the ring of matrices with elements from the ring $\mathbf{K}$. One of the ways to determine determinant of a matrix of $\mathrm{M}(n,\mathbf{K})$ is following ([1, 4, 5]).

**Definition 3.1.** Let the functional

$$\mathrm{d} : \mathrm{M}(n,\mathbf{K}) \to \mathbf{K}$$

satisfy the following axioms.

**Axiom 1.** $\mathrm{d}(\mathbf{A}) = 0$ iff $\mathbf{A}$ is singular (irreversible).

**Axiom 2.** $\forall \mathbf{A}, \mathbf{B} \in \mathrm{M}(n,\mathbf{K})$, $\mathrm{d}(\mathbf{A}\cdot\mathbf{B}) = \mathrm{d}(\mathbf{A})\cdot\mathrm{d}(\mathbf{B})$.

**Axiom 3.** If we obtain a matrix $\mathbf{A}'$ from matrix $\mathbf{A}$ either by adding of an arbitrary row multiplied on the left with its another row or by adding of an arbitrary column multiplied on the right with its another column, then

$$\mathrm{d}(\mathbf{A}') = \mathrm{d}(\mathbf{A})$$

Then the value of the functional d is called determinant of matrix $\mathbf{A} \in \mathrm{M}(n,\mathbf{K})$. □

Known determinants of Dieudonné and Study are examples of such functional. In [1], Aslaksen proved that determinants that satisfy axioms 1, 2, 3, take their value in some commutative subset of the ring. Also it makes no sense for them such property of conventional determinants as the expansion relative to minors of a row or column. Therefore, the determinant representation of the inverse matrix using only these determinants is impossible. This is the reason that causes to introduce determinant functionals that do not satisfy all the above axioms. However in [5], Dyson considers the axiom 1 as necessary to determine the determinant.

In another approach, the determinant of a square matrix over noncommutative ring is considered as a rational function of the entries of matrix. The greatest success is achieved by Gelfand and Retakh in the theory of quasideterminants ([9, 10]). We present introduction into the theory of quasideterminants in section 6.

In third approach, the determinant of a square matrix over noncommutative ring is considered as alternating sum of $n!$ products of the entries of matrix; however, it assumed certain fixed order of factors in each term. E. H. Moore was first who achieved implementation of the key axiom 1 using such definition of noncommutative determinant. Moore had done this not for all square matrices, but only for Hermitian. He defined the determinant of Hermitian matrix[1] $\mathbf{A} = (a_{ij})_{n\times n}$ over division ring with involution by induction over $n$ following way (see [5])

$$\mathrm{Mdet}\mathbf{A} = \left\{ \begin{array}{ll} a_{11}, & n=1 \\ \sum\limits_{j=1}^{n} \varepsilon_{ij} a_{ij} \mathrm{Mdet}\left(\mathbf{A}(i \to j)\right), & n>1 \end{array} \right. \tag{3.1}$$

Here $\varepsilon_{kj} = \left\{ \begin{array}{ll} 1, & i=j \\ -1, & i \neq j \end{array} \right.$, and using $\mathbf{A}(i \to j)$ we denoted the matrix obtained from the matrix $\mathbf{A}$ by successive replacement of column with number $j$ by column with number $i$ and removing from matrix of row with number $i$ and column with number $i$. Another definition of this determinant is presented in [1] using permutations:

$$\mathrm{Mdet}\,\mathbf{A} = \sum_{\sigma\in S_n} |\sigma| a_{n_{11}n_{12}} \cdot \ldots \cdot a_{n_{1l_1}n_{11}} \cdot a_{n_{21}n_{22}} \cdot \ldots \cdot a_{n_{rl_1}n_{r1}}.$$

Here $S_n$ is symmetric group of $n$ elements. A cycle decomposition of a permutation $\sigma$ has form:

$$\sigma = (n_{11}\ldots n_{1l_1})(n_{21}\ldots n_{2l_2})\ldots(n_{r1}\ldots n_{rl_r}).$$

However, there was no any generalization of the definition of determinant by Moore for arbitrary square matrix. Freeman J. Dyson pointed out the importance of this problem in the paper [5].

In papers [2, 3], L. Chen offered the following definition of determinant of square matrix over division ring of quaternions $\mathbf{A} = (a_{ij}) \in \mathrm{M}(n, \mathbf{H})$:

$$\det\mathbf{A} = \sum_{\sigma\in S_n} \varepsilon(\sigma) a_{n_1 i_2} \cdot a_{i_2 i_3} \ldots \cdot a_{i_s n_1} \cdot \ldots \cdot a_{n_r k_2} \cdot \ldots \cdot a_{k_l n_r},$$
$$\sigma = (n_1 i_2 \ldots i_s)\ldots(n_r k_2 \ldots k_r),$$
$$n_1 > i_2, i_3, \ldots, i_s; \ldots, n_r > k_2, k_3, \ldots, k_l,$$
$$n = n_1 > n_2 > \ldots > n_r \geq 1.$$

Despite the fact that this determinant does not satisfy the axiom 1, L. Chen received the determinant representation of the inverse matrix. However we cannot also expand his determinant relative to minors of rows or column (except for $n$-th row). Therefore, L. Chen also did not receive a classical adjoint matrix. For $\mathbf{A} = (\alpha_1, \ldots, \alpha_m)$ over division ring of quaternions $\mathbb{H}$, if $\|\mathbf{A}\| := \det(\mathbf{A}^*\mathbf{A}) \neq 0$, then $\exists \mathbf{A}^{-1} = (b_{jk})$, where

$$\overline{b_{jk}} = \frac{1}{\|\mathbf{A}\|}\omega_{kj}, \quad (j, k = 1, 2, \ldots, n),$$

$$\omega_{kj} = \det(\alpha_1 \ldots \alpha_{j-1}\alpha_n\alpha_{j+1}\ldots\alpha_{n-1}\delta_k)^*(\alpha_1\ldots\alpha_{j-1}\alpha_n\alpha_{j+1}\ldots\alpha_{n-1}\alpha_j).$$

[1]**Hermitian matrix** is such matrix $\mathbf{A} = (a_{ij})$ that $a_{ij} = \overline{a_{ji}}$.

Here $\alpha_i$ is $i$-th column of $\mathbf{A}$, $\delta_k$ is $n$-dimensional column with 1 in $k$-th row and 0 in other ones. L. Chen defined $\|\mathbf{A}\| := \det(\mathbf{A}^*\mathbf{A})$, as double determinant. The solution of the right system of linear equations

$$\sum_{j=1}^{n} \alpha_j x_j = \beta$$

over $\mathbb{H}$ in case $\|\mathbf{A}\| \neq 0$ represented by the following formula, which the author calls Cramer's rule

$$x_j = \|\mathbf{A}\|^{-1}\overline{\mathbf{D}_j},\ \ \forall j = \overline{1,n}$$

where

$$\mathbf{D}_j = \det \begin{pmatrix} \alpha_1^* \\ \vdots \\ \alpha_{j-1}^* \\ \alpha_n^* \\ \alpha_{j+1}^* \\ \vdots \\ \alpha_{n-1}^* \\ \beta^* \end{pmatrix} \begin{pmatrix} \alpha_1 & \dots & \alpha_{j-1} & \alpha_n & \alpha_{j+1} & \dots & \alpha_{n-1} & \alpha_j \end{pmatrix}.$$

Here $\alpha_i$ is $i$-th column of matrix $\mathbf{A}$, $\alpha_i^*$ is $i$-th row of matrix $\mathbf{A}^*$, $\beta^*$ is $n$-dimensional vector-row conjugated with $\beta$.

In this paper, we explore the theory of row and column determinants which develops the classical approach to the definition of the determinant of the matrix, but with a predetermined order of factors in each of the terms of the determinant.

## 4. Quaternion Algebra

A **quaternion algebra** $\mathbb{H}(a,b)$ (we also use notation $\left(\frac{a,b}{\mathbb{F}}\right)$ ) is four-dimensional vector space over the field $\mathbb{F}$ with basis $\{1,i,j,k\}$ and the following multiplication rules:

$$i^2 = a,$$
$$j^2 = b,$$
$$ij = k,$$
$$ji = -k.$$

The field $\mathbb{F}$ is the center of the quaternion algebra $\mathbb{H}(a,b)$.

In the algebra $\mathbb{H}(a,b)$ there are following mappings.

- A quadratic form
  $$\mathrm{n} : x \in \mathbb{H} \to \mathrm{n}(x) \in \mathbb{F}$$
  such that
  $$\mathrm{n}(x \cdot y) = \mathrm{n}(x)\mathrm{n}(y) \quad x, y \in \mathbb{H}$$
  is called the **norm** on quaternion algebra $\mathbb{H}$.
- The linear mapping
  $$\mathrm{t} : x = x^0 + x^1 i + x^2 j + x^3 k \in \mathbb{H} \to \mathrm{t}(x) = 2x^0 \in \mathbb{F}$$

is called **trace of quaternion**. Trace satisfies **permutability property of the trace**:

$$\mathrm{t}\,(q\cdot p)=\mathrm{t}\,(p\cdot q)$$

From the theorem [17]-10.3.3, it follows

$$\mathrm{t}(x)=\frac{1}{2}(x-ixi-jxj-kxk) \tag{4.1}$$

- A linear mapping

$$x\to\overline{x}=\mathrm{t}(x)-x \tag{4.2}$$

is an **involution**. Involution has following properties

$$\overline{\overline{x}}=x$$

$$\overline{x+y}=\overline{x}+\overline{y}$$

$$\overline{x\cdot y}=\overline{y}\cdot\overline{x}$$

A quaternion $\overline{x}$ is called the **conjugate of quaternion** $x\in\mathbb{H}$. Norm and involution satisfy the following condition:

$$\mathrm{n}\,(\overline{q})=\mathrm{n}(q)$$

Trace and involution satisfy the following condition:

$$\mathrm{t}\,(\overline{x})=\mathrm{t}(x)$$

From equations (4.1), (4.2), it follows that

$$\overline{x}=-\frac{1}{2}(x+ixi+jxj+kxk)$$

Depending on the choice of the field $\mathbb{F}$, $a$ and $b$, on the set of quaternion algebra there are only two possibilities [19]:

1. $\left(\frac{a,b}{\mathbb{F}}\right)$ is a division algebra.
2. $\left(\frac{a,b}{\mathbb{F}}\right)$ is isomorphic to the algebra of all $2\times 2$ matrices with entries from the field $\mathbb{F}$. In this case, quaternion algebra is splittable.

Consider some non-isomorphic quaternion algebra with division.

1. $\left(\frac{a,b}{\mathbb{R}}\right)$, where $\mathbb{R}$ is real field, is isomorphic to the Hamilton quaternion skew field **H** whenever $a<0$ and $b<0$. Otherwise $\left(\frac{a,b}{\mathbb{R}}\right)$ is splittable.
2. If $\mathbb{F}$ is the rational field $\mathbb{Q}$, then there exist infinitely many nonisomorphic division quaternion algebras $\left(\frac{a,b}{\mathbb{Q}}\right)$ depending on choice of $a<0$ and $b<0$.
3. Let $\mathbb{Q}_p$ be the $p$-adic field where $p$ is a prime number. For each prime number $p$ there is a unique division quaternion algebra.

## 5. Introduction to the Theory of the Row and Column Determinants

The row and column determinants are introduced for matrices with entries in quaternion algebra $\mathbb{H}$. To introduce the row and column determinants, we need the following definitions in the theory of permutations.

**Definition 5.1.** Let $S_n$ be symmetric group on the set $I_n = \{1, 2, ..., n\}$ ([6], p. 30). If two-line notation of a permutation $\sigma \in S_n$ corresponds to its some cycle notation, then we say that the permutation $\sigma \in S_n$ forms the direct product of disjoint cycles

$$\sigma = \begin{pmatrix} n_{11} & n_{12} & \dots & n_{1l_1} & \dots & n_{r1} & n_{r2} & \dots & n_{rl_r} \\ n_{12} & n_{13} & \dots & n_{11} & \dots & n_{r2} & n_{r3} & \dots & n_{r1} \end{pmatrix} \tag{5.1}$$

□

**Definition 5.2.** We say that cycle notation (5.1) of the permutation $\sigma$ is **left-ordered**, if elements, which closes each of its independent cycles, are recorded first on the left in each cycle

$$\sigma = (n_{11}n_{12} \dots n_{1l_1})\ (n_{21}n_{22} \dots n_{2l_2}) \dots (n_{r1}n_{r2} \dots n_{rl_r})\,.$$

□

**Definition 5.3.** We say that cycle notation (5.1) of the permutation $\sigma$ is **right-ordered**, if elements, which closes each of its independent cycles, are recorded first on the right in each cycle

$$\sigma = (n_{12} \dots n_{1l_1}n_{11})\ (n_{22} \dots n_{2l_2}n_{21}) \dots (n_{r2} \dots n_{rl_r}n_{r1})\,.$$

□

**Definition 5.4.** Let $S_n$ be symmetric group on the set $I_n = \{1, 2, ..., n\}$. Let $i \in \overline{1, n}$. The $i$th **row determinant** of matrix $\mathbf{A} = (a_{ij}) \in \mathrm{M}\,(n, \mathbb{H})$ is defined as expression

$$\mathrm{rdet}_i\mathbf{A} = \sum_{\sigma \in S_n} (-1)^{n-r} a_{i\,i_{k_1}} a_{i_{k_1} i_{k_1+1}} \dots a_{i_{k_1+l_1} i} \dots a_{i_{k_r} i_{k_r+1}} \dots a_{i_{k_r+l_r} i_{k_r}}$$

This expression is defined as the alternative sum of $n!$ monomials compounded from entries of matrix $\mathbf{A}$ such that the index permutation $\sigma \in S_n$ forms the direct product of disjoint independent cycles

$$\sigma = (i\,i_{k_1} i_{k_1+1} \dots i_{k_1+l_1})\,(i_{k_2} i_{k_2+1} \dots i_{k_2+l_2}) \dots (i_{k_r} i_{k_r+1} \dots i_{k_r+l_r})$$

In left-ordered cycle notation of the permutation $\sigma$, the index $i$ starts the first cycle from the left and other cycles satisfy the following conditions

$$i_{k_2} < i_{k_3} < \dots < i_{k_r}, \quad i_{k_t} < i_{k_t+s}, \quad \left(\forall t = \overline{2, r}\right), \quad \left(\forall s = \overline{1, l_t}\right).$$

□

**Definition 5.5.** Let $S_n$ be symmetric group on the set $J_n = \{1, 2, ..., n\}$. Let $j \in \overline{1, n}$. The $j$th **column determinant** of matrix $\mathbf{A} = (a_{ij}) \in \mathrm{M}\,(n, \mathbb{H})$ is defined as expression

$$\mathrm{cdet}_j\,\mathbf{A} = \sum_{\tau \in S_n} (-1)^{n-r}\, a_{j_{k_r} j_{k_r+l_r}} \dots a_{j_{k_r+1} i_{k_r}} \dots a_{j\,j_{k_1+l_1}} \dots a_{j_{k_1+1} j_{k_1}} a_{j_{k_1} j}\,.$$

This expression is defined as the alternative sum of $n!$ monomials compounded from entries of matrix $\mathbf{A}$ such that the index permutation $\tau \in S_n$ forms the direct product of disjoint independent cycles

$$\tau = (j_{k_r+l_r} \dots j_{k_r+1} j_{k_r}) \dots (j_{k_2+l_2} \dots j_{k_2+1} j_{k_2})\ (j_{k_1+l_1} \dots j_{k_1+1} j_{k_1} j)$$

In right-ordered cycle notation of the permutation $\tau$, the index $j$ starts the first cycle from the right and other cycles satisfy the following conditions

$$j_{k_2} < j_{k_3} < \ldots < j_{k_r}, \quad j_{k_t} < j_{k_t+s}, \quad \left(\forall t = \overline{2,r}\right), \quad \left(\forall s = \overline{1,l_t}\right).$$

□

*Remark* 5.6. The peculiarity of the column determinant is that, at the direct calculation, factors of each of the monomials are written from right to left. □

In lemmas 5.7 and 5.8, we consider recursive definition of column and row determinants. This definition is an analogue of the expansion of determinant by rows and columns in commutative case.

*Lemma* 5.7. Let $R_{i\,j}$ be the **right** $ij$**th cofactor** of entry $a_{i\,j}$ of matrix $\mathbf{A} = (a_{ij}) \in \mathrm{M}\,(n,\mathbb{H})$, namely

$$\mathrm{rdet_i}\,\mathbf{A} = \sum_{\mathrm{j}=1}^{\mathrm{n}} \mathrm{a_{i\,j}} \cdot \mathrm{R_{i\,j}} \quad i = \overline{1,n}$$

Then

$$R_{i\,j} = \begin{cases} -\mathrm{rdet}_j\ (\mathbf{A}^{ii}_{.j}(\mathbf{a}_{.\,i})) & i \neq j \\ \mathrm{rdet}_k\ \mathbf{A}^{ii} & i = j \end{cases}$$

$$k = \begin{cases} 2 & i = 1 \\ 1 & i > 1 \end{cases}$$

where we obtain the matrix $(\mathbf{A}^{ii}_{.j}(\mathbf{a}_{.\,i}))$ from $\mathbf{A}$ by replacing its $j$th column with the $i$th column $\mathbf{b}$ and after that by deleting both the $i$th row and the $i$th column. □

*Lemma* 5.8. Let $L_{i\,j}$ be the **left** $ij$**th cofactor** of entry $a_{i\,j}$ of matrix $\mathbf{A} = (a_{ij}) \in \mathrm{M}\,(n,\mathbb{H})$, namely

$$\mathrm{cdet_j}\,\mathbf{A} = \sum_{\mathrm{i}=1}^{\mathrm{n}} \mathrm{L_{i\,j}} \cdot \mathrm{a_{i\,j}} \quad j = \overline{1,n}$$

Then

$$L_{i\,j} = \begin{cases} -\mathrm{cdet}_i\ (\mathbf{A}^{jj}_{i.}(\mathbf{a}_{j\,.})) & i \neq j \\ \mathrm{cdet}_k\ \mathbf{A}^{jj} & i = j \end{cases}$$

$$k = \begin{cases} 2 & j = 1 \\ 1 & j > 1 \end{cases}$$

where we obtain the matrix $(\mathbf{A}^{jj}_{i.}(\mathbf{a}_{j\,.}))$ from $\mathbf{A}$ by replacing its $i$th row with the $j$th and after that by deleting both the $j$th row and the $j$th column. □

*Remark* 5.9. Clearly, any monomial of each row or column determinant of a square matrix corresponds to a certain monomial of another row or column determinant such that both of them have the same sign, consists of the same factors and differ only in their ordering. If the entries of an arbitrary matrix $\mathbf{A}$ are commutative, then $\mathrm{rdet}_1\,\mathbf{A} = \ldots = \mathrm{rdet}_n\mathbf{A} = \mathrm{cdet}_1\,\mathbf{A} = \ldots = \mathrm{cdet}_n\mathbf{A}$. □

Consider the basic properties of the column and row determinants of a square matrix over $\mathbb{H}$. Their proofs immediately follow from the definitions.

**Theorem 5.10.** *If one of the rows (columns) of the matrix* $\mathbf{A} \in \mathrm{M}(n, \mathbb{H})$ *consists of zeros only, then for all* $i = \overline{1, n}$,

$$\mathrm{rdet}_i\, \mathbf{A} = 0, \quad \mathrm{cdet}_i\, \mathbf{A} = 0.$$

**Theorem 5.11.** *If the ith row of the matrix* $\mathbf{A} \in \mathrm{M}(n, \mathbb{H})$ *is left-multiplied by* $b \in \mathbb{H}$, *then for all* $i = \overline{1, n}$

$$\mathrm{rdet}_i\, \mathbf{A}_{i\,.}(b \cdot \mathbf{a}_{.\,i}) = b \cdot \mathrm{rdet}_i\, \mathbf{A}.$$

**Theorem 5.12.** *If the jth column of the matrix* $\mathbf{A} \in \mathrm{M}(n, \mathbb{H})$ *is right-multiplied by* $b \in \mathbb{H}$, *then for all* $j = \overline{1, n}$

$$\mathrm{cdet}_j\, \mathbf{A}_{.\,j}(\mathbf{a}_{j\,.} \cdot b) = \mathrm{cdet}_j\, \mathbf{A} \cdot b.$$

**Theorem 5.13.** *If for* $\mathbf{A} \in \mathrm{M}(n, \mathbb{H})$ *there exists index* $k \in I_n$ *such that* $a_{kj} = b_j + c_j$ *for all* $j = \overline{1, n}$, *then for any* $i = \overline{1, n}$

$$\begin{aligned} \mathrm{rdet}_i\, \mathbf{A} &= \mathrm{rdet}_i\, \mathbf{A}_{k\,.}(\mathbf{b}) + \mathrm{rdet}_i\, \mathbf{A}_{k\,.}(\mathbf{c}), \\ \mathrm{cdet}_i\, \mathbf{A} &= \mathrm{cdet}_i\, \mathbf{A}_{k\,.}(\mathbf{b}) + \mathrm{cdet}_i\, \mathbf{A}_{k\,.}(\mathbf{c}). \end{aligned}$$

*where* $\mathbf{b} = (b_1, \dots, b_n)$, $\mathbf{c} = (c_1, \dots, c_n)$.

**Theorem 5.14.** *If for* $\mathbf{A} \in \mathrm{M}(n, \mathbb{H})$ *there exists index* $k \in I_n$ *such that* $a_{i\,k} = b_i + c_i$ *for all* $i = \overline{1, n}$, *then for any* $j = \overline{1, n}$

$$\begin{aligned} \mathrm{rdet}_j\, \mathbf{A} &= \mathrm{rdet}_j\, \mathbf{A}_{.\,k}(\mathbf{b}) + \mathrm{rdet}_j\, \mathbf{A}_{.\,k}(\mathbf{c}), \\ \mathrm{cdet}_j\, \mathbf{A} &= \mathrm{cdet}_j\, \mathbf{A}_{.\,k}(\mathbf{b}) + \mathrm{cdet}_j \mathbf{A}_{.\,k}(\mathbf{c}). \end{aligned}$$

*where* $\mathbf{b} = (b_1, \dots, b_n)^T$, $\mathbf{c} = (c_1, \dots, c_n)^T$.

**Theorem 5.15.** *Let* $\mathbf{A}^*$ *be the Hermitian adjoint matrix of* $\mathbf{A} \in \mathrm{M}(n, \mathbb{H})$, *then* $\mathrm{rdet}_i\, \mathbf{A}^* = \overline{\mathrm{cdet}_i\, \mathbf{A}}$ *for any* $i = \overline{1, n}$. □

The following theorem is crucial in the theory of the row and column determinants.

**Theorem 5.16.** *If* $\mathbf{A} = (a_{ij}) \in \mathrm{M}(n, \mathbb{H})$ *is a Hermitian matrix, then*

$$\mathrm{rdet}_1 \mathbf{A} = \dots = \mathrm{rdet}_n \mathbf{A} = \mathrm{cdet}_1 \mathbf{A} = \dots = \mathrm{cdet}_n \mathbf{A} \in \mathbb{F}.$$

□

*Remark* 5.17. Since all column and row determinants of a Hermitian matrix over division ring $\mathbb{H}$ are equal, we can define the determinant of a Hermitian matrix $\mathbf{A} \in \mathrm{M}(n, \mathbb{H})$

$$\det \mathbf{A} := \mathrm{rdet}_i\, \mathbf{A} = \mathrm{cdet}_i\, \mathbf{A}, \ \left(\forall i = \overline{1, n}\right).$$

□

*Remark* 5.18. Since we represent determinant of the Hermitian matrix as $i$-th row determinant ($i$ is arbitrary), then according to lemma 5.7, we have

$$\begin{aligned} \det \mathbf{A} &= -\sum_{\sigma \in S_n} a_{i\,j} \cdot \mathrm{rdet}_j\, (\mathbf{A}^{ii}_{.j}(\mathbf{a}_{.\,i})) + a_{i\,i} \cdot \mathrm{rdet}_k\, \mathbf{A}^{ii} \\ k &= \begin{cases} 2 & i = 1 \\ 1 & i > 1 \end{cases} \end{aligned} \tag{5.2}$$

By comparing expressions (3.1) and (5.2) for Hermitian matrix $\mathbf{A} \in \mathrm{M}(n, \mathbf{H})$, we conclude that the row determinant of a Hermitian matrix coincides with the Moore determinant. Hence the row and column determinants extend the Moore determinant to an arbitrary square matrix. □

Properties of the determinant of a Hermitian matrix is completely explored in [14] by its row and column determinants. Among all, consider the following.

**Theorem 5.19.** *If the $i$th row of the Hermitian matrix $\mathbf{A} \in \mathrm{M}(n, \mathbb{H})$ is replaced with a left linear combination of its other rows*

$$\mathbf{a}_{i\,.} = c_1\mathbf{a}_{i_1\,.} + \ldots + c_k\mathbf{a}_{i_k\,.}$$

*where $c_l \in \mathbb{H}$ for all $l = \overline{1,k}$ and $\{i, i_l\} \subset I_n$, then for all $i = \overline{1,n}$*

$$\mathrm{cdet}_i\mathbf{A}_{i\,.}\left(c_1 \cdot \mathbf{a}_{i_1\,.} + \ldots + c_k \cdot \mathbf{a}_{i_k\,.}\right) = \mathrm{rdet}_i\mathbf{A}_{i\,.}\left(c_1 \cdot \mathbf{a}_{i_1\,.} + \ldots + c_k \cdot \mathbf{a}_{i_k\,.}\right) = 0.$$

□

**Theorem 5.20.** *If the $j$th column of a Hermitian matrix $\mathbf{A} \in \mathrm{M}(n, \mathbb{H})$ is replaced with a right linear combination of its other columns*

$$\mathbf{a}_{.\,j} = \mathbf{a}_{.\,j_1}c_1 + \ldots + \mathbf{a}_{.\,j_k}c_k$$

*where $c_l \in \mathbb{H}$ for all $l = \overline{1,k}$ and $\{j, j_l\} \subset J_n$ , then for all $j = \overline{1,n}$*

$$\mathrm{cdet}_j\mathbf{A}_{.\,j}\left(\mathbf{a}_{.\,j_1} \cdot c_1 + \ldots + \mathbf{a}_{.\,j_k} \cdot c_k\right) = \mathrm{rdet}_j\mathbf{A}_{.\,j}\left(\mathbf{a}_{.\,j_1} \cdot c_1 + \ldots + \mathbf{a}_{.\,j_k} \cdot c_k\right) = 0.$$

□

The following theorem on the determinantal representation of the matrix inverse of the Hermitian one follows directly from these properties.

**Theorem 5.21.** *There exist a unique right inverse matrix $(R\mathbf{A})^{-1}$ and a unique left inverse matrix $(L\mathbf{A})^{-1}$ of a nonsingular Hermitian matrix $\mathbf{A} \in \mathrm{M}(n, \mathbb{H})$, ($\det\mathbf{A} \neq 0$), where $(R\mathbf{A})^{-1} = (L\mathbf{A})^{-1} =: \mathbf{A}^{-1}$. Right inverse and left inverse matrices has following determinantal representation*

$$(R\mathbf{A})^{-1} = \frac{1}{\det\mathbf{A}}\begin{pmatrix} R_{11} & R_{21} & \cdots & R_{n1} \\ R_{12} & R_{22} & \cdots & R_{n2} \\ \cdots & \cdots & \cdots & \cdots \\ R_{1n} & R_{2n} & \cdots & R_{nn} \end{pmatrix}$$

$$(L\mathbf{A})^{-1} = \frac{1}{\det\mathbf{A}}\begin{pmatrix} L_{11} & L_{21} & \cdots & L_{n1} \\ L_{12} & L_{22} & \cdots & L_{n2} \\ \cdots & \cdots & \cdots & \cdots \\ L_{1n} & L_{2n} & \cdots & L_{nn} \end{pmatrix}.$$

*where $R_{ij}$, $L_{ij}$ are right and left $ij$-th cofactors of $\mathbf{A}$ respectively for all $i, j = \overline{1,n}$.* □

To obtain the determinantal representation for an arbitrary inverse matrix over division ring $\mathbb{H}$, we consider the right $\mathbf{A}\mathbf{A}^*$ and left $\mathbf{A}^*\mathbf{A}$ corresponding Hermitian matrix.

**Theorem 5.22** ([14])**.** *If an arbitrary column of matrix $\mathbf{A} \in \mathbb{H}^{m\times n}$ is a right linear combination of its other columns, or an arbitrary row of matrix $\mathbf{A}^*$ is a left linear combination of its other rows, then $\det\mathbf{A}^*\mathbf{A} = 0$.* □

Since the principal submatrices of a Hermitian matrix are also Hermitian, then the basis principal minor may be defined in this noncommutative case as a principal nonzero minor of maximal order. We also introduce a notion the **rank of a Hermitian matrix by principal minors**, as maximal order of a principal nonzero minor. The following theorem establishes the correspondence between the rank by principal minors of a Hermitian matrix and the rank of matrix defined as the maximum number of right-linearly independent columns or left-linearly independent rows, which form a basis.

**Theorem 5.23** ([14])**.** *A rank by principal minors of a Hermitian matrix* $\mathbf{A}^*\mathbf{A}$ *is equal to its rank and a rank of* $\mathbf{A} \in \mathbb{H}^{m\times n}$*.* □

**Theorem 5.24** ([14])**.** *If* $\mathbf{A} \in \mathbb{H}^{m\times n}$*, then arbitrary column of matrix* $\mathbf{A}$ *is a right linear combination of its basic columns or arbitrary row of matrix* $\mathbf{A}$ *is a left linear combination of its basic rows.* □

This implies a criterion for the singularity of a corresponding Hermitian matrix.

**Theorem 5.25** ([14])**.** *The right linearly independence of columns of the matrix* $\mathbf{A} \in \mathbb{H}^{m\times n}$ *or the left linearly independence of rows of* $\mathbf{A}^*$ *is the necessary and sufficient condition for*

$$\det \mathbf{A}^*\mathbf{A} \neq 0$$

**Theorem 5.26** ([14])**.** *If* $\mathbf{A} \in \mathrm{M}(n, \mathbb{H})$*, then* $\det \mathbf{A}\mathbf{A}^* = \det \mathbf{A}^*\mathbf{A}$*.* □

**Example 5.27.** Consider the matrix

$$\mathbf{A} = \begin{pmatrix} a_{11} & a_{12} \\ a_{21} & a_{22} \end{pmatrix}$$

Then

$$\mathbf{A}^* = \begin{pmatrix} \overline{a_{11}} & \overline{a_{21}} \\ \overline{a_{12}} & \overline{a_{22}} \end{pmatrix}$$

Respectively, we have

$$\mathbf{A}\mathbf{A}^* = \begin{pmatrix} a_{11}\overline{a_{11}} + a_{12}\overline{a_{12}} & a_{11}\overline{a_{21}} + a_{12}\overline{a_{22}} \\ a_{21}\overline{a_{11}} + a_{22}\overline{a_{12}} & a_{21}\overline{a_{21}} + a_{22}\overline{a_{22}} \end{pmatrix}$$

$$\mathbf{A}^*\mathbf{A} = \begin{pmatrix} \overline{a_{11}}a_{11} + \overline{a_{21}}a_{21} & \overline{a_{11}}a_{12} + \overline{a_{21}}a_{22} \\ \overline{a_{12}}a_{11} + \overline{a_{22}}a_{21} & \overline{a_{12}}a_{12} + \overline{a_{22}}a_{22} \end{pmatrix}$$

To assess the determinants of obtained Hermitian matrices, we consider for example the first row determinant of each matrix. According to the theorem 5.16 and the remark 5.17 we have

$$\det \mathbf{A}\mathbf{A}^* = \mathrm{rdet}_1 \mathbf{A}\mathbf{A}^*$$
$$\det \mathbf{A}^*\mathbf{A} = \mathrm{rdet}_1 \mathbf{A}^*\mathbf{A}$$

According to the lemma 5.7

$$
\begin{aligned}
\det \mathbf{A}\mathbf{A}^* &= (\mathbf{A}\mathbf{A}^*)_{11}(\mathbf{A}\mathbf{A}^*)_{22} - (\mathbf{A}\mathbf{A}^*)_{12}(\mathbf{A}\mathbf{A}^*)_{21}\\
&= (a_{11}\overline{a_{11}} + a_{12}\overline{a_{12}})(a_{21}\overline{a_{21}} + a_{22}\overline{a_{22}})\\
&\quad -(a_{11}\overline{a_{21}} + a_{12}\overline{a_{22}})(a_{21}\overline{a_{11}} + a_{22}\overline{a_{12}})\\
&= a_{11}\overline{a_{11}}a_{21}\overline{a_{21}} + a_{12}\overline{a_{12}}a_{21}\overline{a_{21}}\\
&\quad +a_{11}\overline{a_{11}}a_{22}\overline{a_{22}} + a_{12}\overline{a_{12}}a_{22}\overline{a_{22}}\\
&\quad -a_{11}\overline{a_{21}}a_{21}\overline{a_{11}} - a_{12}\overline{a_{22}}a_{21}\overline{a_{11}}\\
&\quad -a_{11}\overline{a_{21}}a_{22}\overline{a_{12}} - a_{12}\overline{a_{22}}a_{22}\overline{a_{12}}\\
&= a_{12}\overline{a_{12}}a_{21}\overline{a_{21}} + a_{11}\overline{a_{11}}a_{22}\overline{a_{22}}\\
&\quad -a_{12}\overline{a_{22}}a_{21}\overline{a_{11}} - a_{11}\overline{a_{21}}a_{22}\overline{a_{12}}
\end{aligned} \tag{5.3}
$$

$$
\begin{aligned}
\det \mathbf{A}^*\mathbf{A} &= (\mathbf{A}^*\mathbf{A})_{11}(\mathbf{A}^*\mathbf{A})_{22} - (\mathbf{A}^*\mathbf{A})_{12}(\mathbf{A}^*\mathbf{A})_{21}\\
&= (\overline{a_{11}}a_{11} + \overline{a_{21}}a_{21})(\overline{a_{12}}a_{12} + \overline{a_{22}}a_{22})\\
&\quad -(\overline{a_{11}}a_{12} + \overline{a_{21}}a_{22})(\overline{a_{12}}a_{11} + \overline{a_{22}}a_{21})\\
&= \overline{a_{11}}a_{11}\overline{a_{12}}a_{12} + \overline{a_{21}}a_{21}\overline{a_{12}}a_{12}\\
&\quad +\overline{a_{11}}a_{11}\overline{a_{22}}a_{22} + \overline{a_{21}}a_{21}\overline{a_{22}}a_{22}\\
&\quad -\overline{a_{11}}a_{12}\overline{a_{12}}a_{11} - \overline{a_{21}}a_{22}\overline{a_{12}}a_{11}\\
&\quad -\overline{a_{11}}a_{12}\overline{a_{22}}a_{21} - \overline{a_{21}}a_{22}\overline{a_{22}}a_{21}\\
&= \overline{a_{21}}a_{21}\overline{a_{12}}a_{12} + \overline{a_{11}}a_{11}\overline{a_{22}}a_{22}\\
&\quad -\overline{a_{21}}a_{22}\overline{a_{12}}a_{11} - \overline{a_{11}}a_{12}\overline{a_{22}}a_{21}
\end{aligned} \tag{5.4}
$$

Positive terms in the equations (5.3), (5.4) are real numbers and they obviously coincide. To prove equation

$$
a_{12}\overline{a_{22}}a_{21}\overline{a_{11}} + a_{11}\overline{a_{21}}a_{22}\overline{a_{12}} = \overline{a_{21}}a_{22}\overline{a_{12}}a_{11} + \overline{a_{11}}a_{12}\overline{a_{22}}a_{21} \tag{5.5}
$$

we use the rearrangement property of the trace of elements of the quaternion algebra, $\mathrm{t}(pq) = \mathrm{t}(qp)$. Indeed,

$$
a_{12}\overline{a_{22}}a_{21}\overline{a_{11}} + a_{11}\overline{a_{21}}a_{22}\overline{a_{12}} = a_{12}\overline{a_{22}}a_{21}\overline{a_{11}} + \overline{a_{12}\overline{a_{22}}a_{21}\overline{a_{11}}} = \mathrm{t}(a_{12}\overline{a_{22}}a_{21}\overline{a_{11}}),
$$

$$
\overline{a_{21}}a_{22}\overline{a_{12}}a_{11} + \overline{a_{11}}a_{12}\overline{a_{22}}a_{21} = \overline{\overline{a_{11}}a_{12}\overline{a_{22}}a_{21}} + \overline{a_{11}}a_{12}\overline{a_{22}}a_{21} = \mathrm{t}(\overline{a_{11}}a_{12}\overline{a_{22}}a_{21})
$$

Then by the rearrangement property of the trace, we obtain (5.5). □

According to the theorem 5.26 we introduce the concept of double determinant. For the first time this concept was introduced by L. Chen ([2]).

**Definition 5.28.** Determinant of the Hermitian matrix $\mathbf{A}\mathbf{A}^*$ is called **double determinant** of the matrix $\mathbf{A} \in \mathrm{M}(n, \mathbb{H})$

$$
\mathrm{ddet}\mathbf{A} := \det(\mathbf{A}^*\mathbf{A}) = \det(\mathbf{A}\mathbf{A}^*)
$$

□

If $\mathbb{H}$ is the classical quaternion skew field $\mathbf{H}$ over the real field, then the following theorem establishes the validity of Axiom 1 for the double determinant.

**Theorem 5.29.** *If* $\{\mathbf{A}, \mathbf{B}\} \subset \mathrm{M}(n, \mathbf{H})$, *then* $\mathrm{ddet}(\mathbf{A} \cdot \mathbf{B}) = \mathrm{ddet}\mathbf{A} \cdot \mathrm{ddet}\mathbf{B}$. □

Unfortunately, if non-Hermitian matrix is not full rank, then nothing can be said about singularity of its row and column determinant. We show it in the following example.

**Example 5.30.** Consider the matrix

$$A=\begin{pmatrix}i & j\\ j & -i\end{pmatrix}$$

Its second row is obtained from the first row by left-multiplying by $k$. Then, by Theorem 5.25, $\mathrm{ddet}A=0$. Indeed,

$$A^*A=\begin{pmatrix}-i & -j\\ -j & i\end{pmatrix}\cdot\begin{pmatrix}i & j\\ j & -i\end{pmatrix}=\begin{pmatrix}2 & -2k\\ 2k & 2\end{pmatrix}.$$

Then $\mathrm{ddet}A=4+4k^2=0$; however

$$\mathrm{cdet}_1A=\mathrm{cdet}_2A=\mathrm{rdet}_1A=\mathrm{rdet}_2A=-i^2-j^2=2$$

At the same time $rankA=1$, that corresponds to the theorem 5.23. □

The correspondence between the double-determinant and the non-commutative determinants of Moore, Stady and Dieudonné are obtained

$$\mathrm{ddet}\mathbf{A}=\mathrm{Mdet}\left(\mathbf{A}^*\mathbf{A}\right)=\mathrm{Sdet}\mathbf{A}=\mathrm{Ddet}^2\mathbf{A}$$

**Definition 5.31.** Let

$$\mathrm{ddet}\mathbf{A}=\mathrm{cdet}_j\left(\mathbf{A}^*\mathbf{A}\right)=\sum_i\mathbb{L}_{ij}\cdot a_{ij}\quad\forall j=\overline{1,n}$$

then $\mathbb{L}_{ij}$ is called the **left double $ij$th cofactor** of entry $a_{ij}$ of the matrix $\mathbf{A}\in\mathrm{M}\left(n,\mathbb{H}\right)$. □

**Definition 5.32.** Let

$$\mathrm{ddet}\mathbf{A}=\mathrm{rdet}_i\left(\mathbf{A}\mathbf{A}^*\right)=\sum_j a_{ij}\cdot\mathbb{R}_{ij}\quad\forall i=\overline{1,n}$$

then $\mathbb{R}_{ij}$ is called the **right double $ij$th cofactor** of entry $a_{ij}$ of the matrix $\mathbf{A}\in\mathrm{M}\left(n,\mathbb{H}\right)$. □

**Theorem 5.33.** *The necessary and sufficient condition of invertibility of matrix $\mathbf{A}=(a_{ij})\in\mathrm{M}(n,\mathbb{H})$ is $\mathrm{ddet}\mathbf{A}\neq0$. Then $\exists\mathbf{A}^{-1}=\left(L\mathbf{A}\right)^{-1}=\left(R\mathbf{A}\right)^{-1}$, where*

$$(5.6)\qquad \left(L\mathbf{A}\right)^{-1}=\left(\mathbf{A}^*\mathbf{A}\right)^{-1}\mathbf{A}^*=\frac{1}{\mathrm{ddet}\mathbf{A}}\begin{pmatrix}\mathbb{L}_{11} & \mathbb{L}_{21} & \dots & \mathbb{L}_{n1}\\ \mathbb{L}_{12} & \mathbb{L}_{22} & \dots & \mathbb{L}_{n2}\\ \dots & \dots & \dots & \dots\\ \mathbb{L}_{1n} & \mathbb{L}_{2n} & \dots & \mathbb{L}_{nn}\end{pmatrix}$$

$$(5.7)\qquad \left(R\mathbf{A}\right)^{-1}=\mathbf{A}^*\left(\mathbf{A}\mathbf{A}^*\right)^{-1}=\frac{1}{\mathrm{ddet}\mathbf{A}^*}\begin{pmatrix}\mathbb{R}_{11} & \mathbb{R}_{21} & \dots & \mathbb{R}_{n1}\\ \mathbb{R}_{12} & \mathbb{R}_{22} & \dots & \mathbb{R}_{n2}\\ \dots & \dots & \dots & \dots\\ \mathbb{R}_{1n} & \mathbb{R}_{2n} & \dots & \mathbb{R}_{nn}\end{pmatrix}$$

*and* $\mathbb{L}_{ij}=\mathrm{cdet}_j(\mathbf{A}^*\mathbf{A})_{.j}\left(\mathbf{A}^*_{.i}\right)$, $\mathbb{R}_{ij}=\mathrm{rdet}_i(\mathbf{A}\mathbf{A}^*)_{i.}\left(\mathbf{A}^*_{j.}\right)$, $\left(\forall i,j=\overline{1,n}\right)$. □

**Corollary 5.34.** *Let* $\mathrm{ddet}\mathbf{A}\neq0$. *Then*

$$\mathbb{L}_{ij}=\mathbb{R}_{ij}$$

□

*Remark* 5.35. In the theorem 5.33, the inverse matrix $\mathbf{A}^{-1}$ of an arbitrary matrix $\mathbf{A}\in \mathrm{M}(n,\mathbb{H})$ under the assumption of $\mathrm{ddet}\mathbf{A}\neq 0$ is represented by the analog of the classical adjoint matrix. If we denote this analog of the adjoint matrix over $\mathbb{H}$ by $\mathrm{Adj}[[\mathbf{A}]]$, then the next formula is valid over $\mathbb{H}$:

$$\mathbf{A}^{-1}=\frac{\mathrm{Adj}[[\mathbf{A}]]}{\mathrm{ddet}\mathbf{A}}.$$

□

An obvious consequence of a determinantal representation of the inverse matrix by the classical adjoint matrix is the Cramer's rule.

**Theorem 5.36.** *Let*

$$\mathbf{A}\cdot\mathbf{x}=\mathbf{y} \tag{5.8}$$

*be a right system of linear equations with a matrix of coefficients* $\mathbf{A}\in \mathrm{M}(n,\mathbb{H})$*, a column of constants* $\mathbf{y}=(y_1,\dots,y_n)^T\in\mathbb{H}^{n\times 1}$*, and a column of unknowns* $\mathbf{x}=(x_1,\dots,x_n)^T$*. If* $\mathrm{ddet}\mathbf{A}\neq 0$*, then the solution to the system of linear equations* (5.8) *has a unique solution that is represented by equation:*

$$x_j=\frac{\mathrm{cdet}_j\,(\mathbf{A}^*\mathbf{A})_{.j}\,(\mathbf{f})}{\mathrm{ddet}\mathbf{A}}\qquad \forall j=\overline{1,n} \tag{5.9}$$

*where* $\mathbf{f}=\mathbf{A}^*\mathbf{y}$. □

**Theorem 5.37.** *Let*

$$\mathbf{x}\cdot\mathbf{A}=\mathbf{y} \tag{5.10}$$

*be a left system of linear equations with a matrix of coefficients* $\mathbf{A}\in \mathrm{M}(n,\mathbb{H})$*, a column of constants* $\mathbf{y}=(y_1,\dots,y_n)\in\mathbb{H}^{1\times n}$ *and a column of unknowns* $\mathbf{x}=(x_1,\dots,x_n)$*. If* $\mathrm{ddet}\mathbf{A}\neq 0$*, then the system of linear equations* (5.10) *has a unique solution that is represented by equation:*

$$x_i=\frac{\mathrm{rdet}_i\,(\mathbf{A}\mathbf{A}^*)_{i.}\,(\mathbf{z})}{\mathrm{ddet}\mathbf{A}}\qquad \forall i=\overline{1,n} \tag{5.11}$$

*where* $\mathbf{z}=\mathbf{y}\mathbf{A}^*$. □

Equations (5.9) and (5.11) are the obvious and natural generalizations of Cramer's rule for systems of linear equations over quaternion algebra. As follows from the theorem 5.21, the closer analog to Cramer's rule can be obtained in the following specific cases.

**Theorem 5.38.** *Let the matrix of coefficients* $\mathbf{A}\in \mathrm{M}(n,\mathbb{H})$ *in the right system of linear equations* (5.8) *over division ring* $\mathbb{H}$ *be Hermitian. Then the system of linear equations* (5.8) *has the unique solution represented by the equation:*

$$x_j=\frac{\mathrm{cdet}_j\mathbf{A}_{.j}\,(\mathbf{y})}{\det\mathbf{A}},\qquad \left(\forall j=\overline{1,n}\right).$$

□

**Theorem 5.39.** *Let the matrix of coefficients* $\mathbf{A}\in \mathrm{M}(n,\mathbb{H})$ *in the left system of linear equations* (5.10) *over division ring* $\mathbb{H}$ *be Hermitian. Then the system of linear equations* (5.10) *has the unique solution represented by the equation:*

$$x_i=\frac{\mathrm{rdet}_i\mathbf{A}_{i.}\,(\mathbf{y})}{\det\mathbf{A}},\qquad \left(\forall i=\overline{1,n}\right).$$

□

In the framework of the row-column determinants also obtained Cramer's rule for the right $\mathbf{AX}=\mathbf{B}$, left $\mathbf{XA}=\mathbf{B}$ and two-sided $\mathbf{AXB}=\mathbf{C}$ matrix equations are obtained ([15]). Determinantal representations of Moore-Penrose inverse matrix and Cramer's rule for the normal solution of the left and right system of linear equations are given as well ([16]).

## 6. Quasideterminant

**Theorem 6.1.** *Suppose matrix*

$$\mathbf{A}=\begin{pmatrix}a_{11} & ... & a_{1n}\\ ... & ... & ...\\ a_{n1} & ... & a_{nn}\end{pmatrix}$$

*has inverse matrix $A^{-1}$.*[2] *Then minor of inverse matrix satisfies the following equation, provided that the inverse matrices exist*

$$((\mathbf{A}^{-1})_{IJ})^{-1}=\mathbf{A}_{JI}-\mathbf{A}_{J\cdot}^{\cdot I}(\mathbf{A}^{JI})^{-1}\mathbf{A}_{\cdot I}^{J\cdot} \tag{6.1}$$

*Proof.* Definition of inverse matrix leads to the system of linear equations

$$\mathbf{A}^{JI}(\mathbf{A}^{-1})_{\cdot J}^{I\cdot}+\mathbf{A}_{\cdot I}^{J\cdot}(\mathbf{A}^{-1})_{IJ}=0 \tag{6.2}$$

$$\mathbf{A}_{J\cdot}^{\cdot I}(\mathbf{A}^{-1})_{\cdot J}^{I\cdot}+\mathbf{A}_{JI}(\mathbf{A}^{-1})_{IJ}=E \tag{6.3}$$

We multiply (6.2) by $\left(\mathbf{A}^{JI}\right)^{-1}$

$$(\mathbf{A}^{-1})_{\cdot J}^{I\cdot}+(\mathbf{A}^{JI})^{-1}\mathbf{A}_{\cdot I}^{J\cdot}(\mathbf{A}^{-1})_{IJ}=0 \tag{6.4}$$

Now we can substitute (6.4) into (6.3)

$$\mathbf{A}_{JI}(\mathbf{A}^{-1})_{IJ}-\mathbf{A}_{J\cdot}^{\cdot I}(\mathbf{A}^{JI})^{-1}\mathbf{A}_{\cdot I}^{J\cdot}(\mathbf{A}^{-1})_{IJ}=E \tag{6.5}$$

(6.1) follows from (6.5). □

**Corollary 6.2.** *Suppose matrix* $\mathbf{A}$ *has inverse matrix. Then elements of inverse matrix satisfy to the equation*

$$((\mathbf{A}^{-1})_{ij})^{-1}=a_{ji}-\mathbf{A}_{j\cdot}^{\cdot i}(\mathbf{A}^{ji})^{-1}\mathbf{A}_{\cdot i}^{j\cdot} \tag{6.6}$$

□

**Example 6.3.** Consider matrix

$$\mathbf{A}=\begin{pmatrix}a_{11} & a_{12}\\ a_{21} & a_{22}\end{pmatrix}$$

According to (6.6)

$$(A^{-1})_{11}=(a_{11}-a_{12}(a_{22})^{-1}\ a_{21})^{-1} \tag{6.7}$$

$$(A^{-1})_{21}=(a_{21}-a_{22}(a_{12})^{-1}\ a_{11})^{-1} \tag{6.8}$$

$$(A^{-1})_{12}=(a_{12}-a_{11}(a_{21})^{-1}\ a_{22})^{-1} \tag{6.9}$$

$$(A^{-1})_{22}=(a_{22}-a_{21}(a_{11})^{-1}\ a_{12})^{-1} \tag{6.10}$$

□

[2] This statement and its proof are based on statement 1.2.1 from [7] (page 8) for matrix over free division ring.

We call matrix

(6.11) $$\mathcal{H}A=((\mathcal{H}A)_{ij})=((a_{ji})^{-1})$$

**Hadamard inverse of matrix** [3]

$$\mathbf{A}=\begin{pmatrix}a_{11}&...&a_{1n}\\...&...&...\\a_{n1}&...&a_{nn}\end{pmatrix}$$

**Definition 6.4.** $(ji)$**-quasideterminant** of the matrix $\mathbf{A}$ of $n\times n$ matrix $A$ is formal expression

(6.12) $$|\mathbf{A}|_{ji}=(\mathcal{H}\mathbf{A}^{-1})_{ji}=((\mathbf{A}^{-1})_{ij})^{-1}$$

We consider $(ji)$-quasideterminant as an element of the matrix $|A|$, which is called **quasideterminant**. □

**Theorem 6.5.** *Expression for $(ji)$-quasideterminant has form*

(6.13) $$|\mathbf{A}|_{ji}=a_{ji}-\mathbf{A}_{j\cdot}^{\cdot i}(\mathbf{A}^{ji})^{-1}\mathbf{A}_{\cdot i}^{j\cdot}$$

(6.14) $$|\mathbf{A}|_{ji}=a_{ji}-\mathbf{A}_{j\cdot}^{\cdot i}\,\mathcal{H}|\mathbf{A}^{ji}|\,\mathbf{A}_{\cdot i}^{j\cdot}$$

*Proof.* The statement follows from (6.6) and (6.12). □

**Theorem 6.6.** *Let*

(6.15) $$A=\begin{pmatrix}1&0\\0&1\end{pmatrix}$$

*Then*

(6.16) $$A^{-1}=\begin{pmatrix}1&0\\0&1\end{pmatrix}$$

*Proof.* It is clear from (6.7) and (6.10) that $(A^{-1})_{11}=1$ and $(A^{-1})_{22}=1$. However expression for $(A^{-1})_{21}$ and $(A^{-1})_{12}$ cannot be defined from (6.8) and (6.9) since $(a_{21}-a_{22}(a_{12})^{-1}\ a_{11})^{-1}=(a_{12}-a_{11}(a_{21})^{-1}\ a_{22})^{-1}=0$. We can transform these expressions. For instance

$$\begin{aligned}(A^{-1})_{21}&=(a_{21}-a_{22}(a_{12})^{-1}\ a_{11})^{-1}\\&=(a_{11}((a_{11})^{-1}\ a_{12}-(a_{21})^{-1}\ a_{22}))^{-1}\\&=((a_{21})^{-1}\ a_{11}(a_{21}(a_{11})^{-1}\ a_{12}-a_{22}))^{-1}\\&=(a_{11}(a_{21}(a_{11})^{-1}\ a_{12}-a_{22}))^{-1}\ a_{21}\end{aligned}$$

It follows immediately that $(A^{-1})_{21}=0$. In the same manner we can find that $(A^{-1})_{12}=0$. This completes the proof of (6.16). □

From the proof of the theorem 6.6 we see that we cannot always use the equation (6.6) to find elements of inverse matrix and we need more transformations to solve this problem. From the theorem [17]-4.6.3, it follows that if

$$\operatorname{rank}\begin{pmatrix}a_{11}&...&a_{1n}\\...&...&...\\a_{n1}&...&a_{nn}\end{pmatrix}\le n-2$$

[3] [8]-page 4

then $|\mathbf{A}|_{ij}$, $i = 1, ..., n$, $j = 1, ..., n$, is not defined. From this, it follows that although quasideterminant is a powerful tool, use of a determinant is a major advantage.

**Theorem 6.7.** *Let matrix* $\mathbf{A}$ *have inverse matrix. Then for any matrices* $\mathbf{B}$ *and* $\mathbf{C}$ *equation*

$$\mathbf{B} = \mathbf{C} \tag{6.17}$$

*follows from the equation*

$$\mathbf{B}\mathbf{A} = \mathbf{C}\mathbf{A} \tag{6.18}$$

*Proof.* Equation (6.17) follows from the equation (6.18) if we multiply both parts of the equation (6.18) over $\mathbf{A}^{-1}$. □

**Theorem 6.8.** *Solution of nonsingular system of linear equations*

$$\mathbf{A}x = b \tag{6.19}$$

*is determined uniquely and can be presented in either form*[4]

$$x = \mathbf{A}^{-1}b \tag{6.20}$$

$$x = \mathcal{H}|\mathbf{A}|\, b \tag{6.21}$$

*Proof.* Multiplying both sides of equation (6.19) from left by $\mathbf{A}^{-1}$ we get (6.20). Using definition 6.4, we get (6.21). Since theorem 6.7 the solution is unique. □

## 7. Correspondence between Row-Column Determinants and Quasideterminants

**Theorem 7.1.** *If* $\mathbf{A} \in \mathrm{M}(n, \mathbb{H})$ *is an invertible matrix, then, for arbitrary* $p, q = 1, ..., n$, *we have the following representation of quasideterminant*

$$|\mathbf{A}|_{pq} = \mathrm{ddet}\mathbf{A}\ (\mathbb{L}_{\mathrm{pq}})^{-1} = \mathrm{ddet}\mathbf{A}\ (\mathrm{cdet_q}(\mathbf{A}^*\mathbf{A})_{.\,\mathrm{q}}(\mathbf{A}^*_{.\mathrm{p}}))^{-1} \tag{7.1}$$

$$= \frac{\mathrm{ddet}\mathbf{A}}{\mathrm{n}(\mathrm{cdet}_q(\mathbf{A}^*\mathbf{A})_{.\,q}(\mathbf{A}^*_{.p}))}\ \overline{\mathrm{cdet}_q(\mathbf{A}^*\mathbf{A})_{.\,q}(\mathbf{A}^*_{.p})} \tag{7.2}$$

$$|\mathbf{A}|_{pq} = \mathrm{ddet}\mathbf{A}\ (\mathbb{R}_{\mathrm{pq}})^{-1} = \mathrm{ddet}\mathbf{A}\ (\mathrm{rdet_p}(\mathbf{A}\mathbf{A}^*)_{\mathrm{p}\,.}(\mathbf{A}^*_{\mathrm{q}\,.}))^{-1} \tag{7.3}$$

$$= \frac{\mathrm{ddet}\mathbf{A}}{\mathrm{n}(\mathrm{rdet}_p(\mathbf{A}\mathbf{A}^*)_{p\,.}(\mathbf{A}^*_{q\,.}))}\ \overline{\mathrm{rdet}_p(\mathbf{A}\mathbf{A}^*)_{p\,.}(\mathbf{A}^*_{q\,.})} \tag{7.4}$$

*Proof.* Let $\mathbf{A}^{-1} = (b_{ij})$ be matrix inversed to the matrix $\mathbf{A}$. The equation (6.12) reveals the relationship between quasideterminant $|\ \mathbf{A}\ |_{p,q}$ of matrix $\mathbf{A} \in \mathrm{M}(n, \mathbb{H})$ and elements of the inverse matrix $\mathbf{A}^{-1} = (b_{ij})$, namely

$$|\ \mathbf{A}\ |_{pq} = b_{qp}^{-1}$$

for every $p$, $q = 1, ..., n$. At the same time, the theory of row and column determinants (the theorem 5.33) gives us representation of inverse matrix through its left (5.6) and right (5.7) double cofactors. Thus, accordingly, we obtain

$$|\ \mathbf{A}\ |_{pq} = b_{qp}^{-1} = \left(\frac{\mathbb{L}_{pq}}{\mathrm{ddet}\mathbf{A}}\right)^{-1} = \left(\frac{\mathrm{cdet}_q(\mathbf{A}^*\mathbf{A})_{.\,q}\,(\mathbf{A}^*_{.p})}{\mathrm{ddet}\mathbf{A}}\right)^{-1}, \tag{7.5}$$

[4] See similar statement in the theorem [7]-1.6.1.

$$|\,\mathbf{A}\,|_{pq} = b_{qp}^{-1} = \left(\frac{\mathbb{R}_{pq}}{\mathrm{ddet}\mathbf{A}}\right)^{-1} = \left(\frac{\mathrm{rdet}_p(\mathbf{A}\mathbf{A}^*)_{p\,.}\left(\mathbf{A}^*_{q\,.}\right)}{\mathrm{ddet}\mathbf{A}}\right)^{-1}. \tag{7.6}$$

Since $\mathrm{ddet}\mathbf{A} \neq 0 \in \mathbb{F}$, then $\exists(\mathrm{ddet}\mathbf{A})^{-1} \in \mathbb{F}$. It follows that

$$\mathrm{cdet}_q(\mathbf{A}^*\mathbf{A})_{.\,q}\left(\mathbf{A}^*_{.p}\right)^{-1} = \frac{\overline{\mathrm{cdet}_q(\mathbf{A}^*\mathbf{A})_{.\,q}\left(\mathbf{A}^*_{.p}\right)}}{\mathrm{n}(\mathrm{cdet}_q(\mathbf{A}^*\mathbf{A})_{.\,q}\left(\mathbf{A}^*_{.p}\right))}, \tag{7.7}$$

$$\mathrm{rdet}_p(\mathbf{A}\mathbf{A}^*)_{p\,.}\left(\mathbf{A}^*_{q\,.}\right)^{-1} = \frac{\overline{\mathrm{rdet}_p(\mathbf{A}\mathbf{A}^*)_{p\,.}\left(\mathbf{A}^*_{q\,.}\right)}}{\mathrm{n}(\mathrm{rdet}_p(\mathbf{A}\mathbf{A}^*)_{p\,.}\left(\mathbf{A}^*_{q\,.}\right))}. \tag{7.8}$$

Substituting (7.7) into (7.5), and (7.8) into (7.6), we accordingly obtain (7.2) and (7.4).

We proved the theorem. □

Equation (7.2) gives an explicit representation of quasideterminant $|\,\mathbf{A}\,|_{p,q}$ of the matrix $\mathbf{A} \in \mathrm{M}(n,\mathbb{H})$ for all $p,q = 1,...,n$ through the column determinant of its corresponding left Hermitian matrix $\mathbf{A}^*\mathbf{A}$, and (7.4) does through the row determinant of its corresponding right Hermitian matrix $\mathbf{A}\mathbf{A}^*$.

**Example 7.2.** Consider a matrix

$$\mathbf{A} = \begin{pmatrix} a_{11} & a_{12} \\ a_{21} & a_{22} \end{pmatrix}$$

According to (6.13)

$$|A| = \begin{pmatrix} a_{11} - a_{12}(a_{22})^{-1}\,a_{21} & a_{12} - a_{11}(a_{21})^{-1}\,a_{22} \\ a_{21} - a_{22}(a_{12})^{-1}\,a_{11} & a_{22} - a_{21}(a_{11})^{-1}\,a_{12} \end{pmatrix} \tag{7.9}$$

Our goal is to find this quasideterminant, using the theorem 7.1. It is evident that

$$\mathbf{A}^* = \begin{pmatrix} \overline{a_{11}} & \overline{a_{21}} \\ \overline{a_{12}} & \overline{a_{22}} \end{pmatrix} \qquad \mathbf{A}^*\mathbf{A} = \begin{pmatrix} \mathrm{n}(a_{11}) + \mathrm{n}(a_{21}) & \overline{a_{11}}a_{12} + \overline{a_{21}}a_{22} \\ \overline{a_{12}}a_{11} + \overline{a_{22}}a_{21} & \mathrm{n}(a_{12}) + \mathrm{n}(a_{22}) \end{pmatrix}.$$

Calculate the necessary determinants

$$\begin{aligned}
\mathrm{ddet}\mathbf{A} &= \mathrm{rdet}_1(\mathbf{A}^*\mathbf{A}) \\
&= (\mathrm{n}(a_{11}) + \mathrm{n}(a_{21}))\cdot(\mathrm{n}(a_{12}) + \mathrm{n}(a_{22})) \\
&\quad -(\overline{a_{11}}a_{12} + \overline{a_{21}}a_{22})\cdot(\overline{a_{12}}a_{11} + \overline{a_{22}}a_{21}) \\
&= \mathrm{n}(a_{11})\mathrm{n}(a_{12}) + \mathrm{n}(a_{11})\mathrm{n}(a_{22}) + \mathrm{n}(a_{21})\mathrm{n}(a_{12}) + \mathrm{n}(a_{21})\mathrm{n}(a_{22}) \\
&\quad -\overline{a_{11}}a_{12}\overline{a_{12}}a_{11} - \overline{a_{11}}a_{12}\overline{a_{22}}a_{21} - \overline{a_{21}}a_{22}\overline{a_{12}}a_{11} - \overline{a_{21}}a_{22}\overline{a_{22}}a_{21} \\
&= \mathrm{n}(a_{11})\mathrm{n}(a_{22}) + \mathrm{n}(a_{21})\mathrm{n}(a_{12}) - (\overline{a_{11}}a_{12}\overline{a_{22}}a_{21} + \overline{\overline{a_{11}}a_{12}\overline{a_{22}}a_{21}}) \\
&= \mathrm{n}(a_{11})\mathrm{n}(a_{22}) + \mathrm{n}(a_{21})\mathrm{n}(a_{12}) - \mathrm{t}(\overline{a_{11}}a_{12}\overline{a_{22}}a_{21})
\end{aligned}$$

$$\begin{aligned}
\mathrm{cdet}_1(\mathbf{A}^*\mathbf{A})_{.1}(\mathbf{a}^*_{.2}) &= \mathrm{cdet}_1\begin{pmatrix} \overline{a_{21}} & \overline{a_{11}}a_{12} + \overline{a_{21}}a_{22} \\ \overline{a_{22}} & \mathrm{n}(a_{12}) + \mathrm{n}(a_{22}) \end{pmatrix} \\
&= \mathrm{n}(a_{12})\overline{a_{21}} + \mathrm{n}(a_{22})\overline{a_{21}} - \overline{a_{11}}a_{12}\overline{a_{22}} - \overline{a_{21}}a_{22}\overline{a_{22}} \\
&= \mathrm{n}(a_{12})\overline{a_{21}} - \overline{a_{11}}a_{12}\overline{a_{22}}.
\end{aligned}$$

Then

$$\overline{\mathrm{cdet}_1(\mathbf{A}^*\mathbf{A})_{.1}(\mathbf{a}^*_{.2})} = \mathrm{n}(a_{12})a_{21} - a_{22}\overline{a_{12}}a_{11},$$

$$\begin{aligned}
\mathrm{n}(\mathrm{cdet}_1(\mathbf{A}^*\mathbf{A})_{.1}(\mathbf{a}^*_{.2})) &= \overline{\mathrm{cdet}_1(\mathbf{A}^*\mathbf{A})_{.1}(\mathbf{a}^*_{.2})}\cdot \mathrm{cdet}_1(\mathbf{A}^*\mathbf{A})_{.1}(\mathbf{a}^*_{.2})\\
&= (\mathrm{n}(a_{12})a_{21} - a_{22}\overline{a_{12}}a_{11})\cdot(\mathrm{n}(a_{12})\overline{a_{21}} - \overline{a_{11}}a_{12}\overline{a_{22}})\\
&= \mathrm{n}^2(a_{12})\mathrm{n}(a_{21}) - \mathrm{n}(a_{12})a_{21}\overline{a_{11}}a_{12}\overline{a_{22}}\\
&\quad -\mathrm{n}(a_{12})a_{22}\overline{a_{12}}a_{11}\overline{a_{21}} + a_{22}\overline{a_{12}}a_{11}\overline{a_{11}}a_{12}\overline{a_{22}}\\
&= \mathrm{n}(a_{12})(\mathrm{n}(a_{12})\mathrm{n}(a_{21}) - \mathrm{t}(\overline{a_{11}}a_{12}\overline{a_{22}}a_{21}) + \mathrm{n}(a_{21})\mathrm{n}(a_{12}))\\
&= \mathrm{n}(a_{12})\mathrm{ddet}\mathbf{A}.
\end{aligned}$$

Following the formula (7.2), we obtain

$$\begin{aligned}
|\mathbf{A}|_{21} &= \frac{\mathrm{ddet}\mathbf{A}}{\mathrm{n}(\mathrm{cdet}_1(\mathbf{A}^*\mathbf{A})_{.1}(\mathbf{a}^*_{.2}))}\overline{\mathrm{cdet}_1(\mathbf{A}^*\mathbf{A})_{.1}(\mathbf{a}^*_{.2})}\\
&= \frac{\mathrm{ddet}\mathbf{A}}{\mathrm{n}(a_{12})\mathrm{ddet}\mathbf{A}}\overline{\mathrm{cdet}_1(\mathbf{A}^*\mathbf{A})_{.1}(\mathbf{a}^*_{.2})}\\
&= \frac{1}{\mathrm{n}(a_{12})}\cdot\overline{\mathrm{cdet}_1(\mathbf{A}^*\mathbf{A})_{.1}(\mathbf{a}^*_{.2})}\\
&= \frac{1}{\mathrm{n}(a_{12})}\cdot(\mathrm{n}(a_{12})a_{21} - a_{22}\overline{a_{12}}a_{11})\\
&= a_{21} - a_{22}(a_{12})^{-1}a_{11}.
\end{aligned} \tag{7.10}$$

The last expression in (7.10) coincides with the expression $|\mathbf{A}|_{21}$ in (7.9). □

## 9. Index

# 10. Special Symbols and Notations

$\mathbf{A}_{.\,T}$ minor 2
$\mathbf{A}_{S\,.}$ minor 2
$\mathbf{A}^{.\,T}$ minor 3
$\mathbf{A}^{S\,.}$ minor 2

$|\mathbf{A}|_{ji}$ $(ji)$-quasideterminant of matrix $\mathbf{A}$ 16
$\mathbf{A}_{.j}\left(\mathbf{b}\right)$ transformation of matrix 3
$\mathbf{A}_{i.}\left(\mathbf{b}\right)$ transformation of matrix 3
$|A|$ quasideterminant of matrix $\mathbf{A}$ 16

$\mathrm{cdet_j}\,\mathbf{A}$ $j$th column determinant of matrix $\mathbf{A}$ 7

$\mathrm{ddet}\mathbf{A}$ double determinant of matrix $\mathbf{A}$ 12

$\left(\frac{a,b}{\mathbb{F}}\right)$ quaternion algebra over field $\mathbb{F}$ 5

$\mathcal{H}A$ Hadamard inverse of matrix 16
$\mathbb{H}(a,b)$ quaternion algebra 5

$L_{i\,j}$ left $ij$th cofactor of entry of matrix 8
$\mathbb{L}_{ij}$ left double $ij$th cofactor of entry of matrix 13

$\mathrm{n}(x)$ norm of quaternion $x$ 5

$R_{i\,j}$ right $ij$th cofactor of entry of matrix 8

$\mathbb{R}_{ij}$ right double $ij$th cofactor of entry of matrix 13
$\mathrm{rdet}_i\mathbf{A}$ $i$th row determinant of matrix $\mathbf{A}$ 7

$S_n$ symmetric group 7

$\mathrm{t}(x)$ trace of quaternion $x$ 5

# Соответствие между строчно-столбцовыми определителями и квазидетерминантами матриц над кватернионной алгеброй.

Александр Клейн and Иван Кирчей

Аннотация. В работе рассматриваются элементы теорий квазидетерминантов и строчно-столбцовых определителей и их применение к решению систем линейных уравнений в кватернионной алгебре. Установлено соответствие между строчно-столбцовыми определителями и квазидетерминантами матриц над кватернионной алгеброй.

## Содержание

## 1. Предисловие

Линейная алгебра является мощным инструментом, которым мы пользуемся в различных областях математики, включая математический анализ, аналитическую и дифференциальную геометрии, теорию дифференциальных уравнений, теорию оптимального управления. Линейная алгебра накопила богатый набор различных методов. Так как некоторые методы имеют общий конечный результат, это даёт нам возможность выбрать наиболее эффективный метод в зависимости от характера расчётов или выполняемых построений.

Aleks_Kleyn@MailAPS.org.

kyrchei@lms.lviv.ua.

Институт Прикладных Проблем Механики и Математики им. Я.С.Пидстрыгача НАН Украины, ул. Наукова 3а, Львов, 79053, Украина.

При переходе от линейной алгебры над полем к линейной алгебре над телом мы хотим сохранить как можно больше инструментария, которым мы регулярно пользуемся. Уже в начале XX века, вскоре после после создания Гамильтоном алгебры кватернионов, математики стали искать ответ, как выглядит алгебра с некоммутативным умножением. В частности возникла проблема как определить детерминант матрицы с элементами, принадлежащими некоммутативному кольцу. Такой детерминант также называют некоммутативным детерминантом.

Возникло немало подходов к определению некоммутативного детерминанта. Но ни один из введённых некоммутативных детерминантов в полной мере не сохранял те свойства, которыми он обладал для матриц над полем. Более того, в работе [11] доказано, что не существует единого детерминантного функционала, который бы расширял определение детерминанта вещественных матриц для матриц над телом кватернионов. Поэтому поиски решения проблемы определения некоммутативного детерминанта продолжаются до сих пор.

В данной работе мы рассматриваем два подхода к определению некоммутативного детерминанта: строчно-столбцовые определители и квазидетерминанты.

Строчно-столбцовые определители являются расширением классического определения детерминанта, но с заранее определенным порядком элементов в каждом из слагаемых детерминанта. При использовании строчно-столбцовых определителей, мы получаем решение системы линейных уравнений над кватернионной алгеброй согласно правилу Крамера.

Квазидетерминант появился в результате анализа процесса обращения матрицы. При использовании квазидетерминанта, решение системы линейных уравнений над кватернионной алгеброй аналогично методу Гаусса.

Общим в определении строчно-столбцовых определителей и квазидетерминантов является, что в обоих случаях квадратной матрице $n$-го порядка с некоммутативными элементами ставится в соответствие не один определитель, а некоторое их множество ($n^2$ квазидетерминантов и соответственно $n$ строчных и $n$ столбцовых определителей).

Квазидетерминанты нашли в последнее время широкое применение, как в области линейной алгебры ([18, 21]), так и физики ([12, 13, 20]). Сравнительно недавно введенные строчно-столбцовые определители ([14]) менее известны. Целью данной статьи является установление соответствия между строчно-столбцовыми определителями и квазидетерминантами матриц над кватернионной алгеброй. Авторы полны надежды что установленное соответствие может дать обоюдное развитие как теории квазидетерминантов, так и строчно-столбцовых определителей.

## 2. Соглашение об обозначениях

Существуют различные формы записи элементов матрицы. В данной статье элемент матрицы $\mathbf{A}$ мы обозначать $a_{ij}$. При этом индекс $i$ нумерует строки, а индекс $j$ нумерует столбцы.

Для обозначения различных миноров матрицы $\mathbf{A}$ мы будем пользоваться следующей записью

$\mathbf{a}_{i\,\cdot}$ **:** строка с индексом $i$

$\mathbf{A}_{S\,\cdot}$ : минор, полученный из $A$ выбором строк с индексом из множества $S$
$\mathbf{A}^{i\,\cdot}$ : минор, полученный из $A$ удалением строки $\mathbf{a}_{i\,\cdot}$
$\mathbf{A}^{S\,\cdot}$ : минор, полученный из $A$ удалением строк с индексом из множества $S$
$\mathbf{a}_{\cdot\,j}$ : столбец с индексом $j$
$\mathbf{A}_{\cdot\,T}$ : минор, полученный из $A$ выбором столбцов с индексом из множества $T$
$\mathbf{A}^{\cdot\,j}$ : минор, полученный из $A$ удалением столбца $\mathbf{a}_{\cdot\,j}$
$\mathbf{A}^{\cdot\,T}$ : минор, полученный из $A$ удалением столбцов с индексом из множества $T$
$\mathbf{A}_{.j}\,(\mathbf{b})$ : матрица, которую получим из матрицы $\mathbf{A}$ заменой ее $j$-го столбца столбцом $\mathbf{b}$
$\mathbf{A}_{i.}\,(\mathbf{b})$ : матрица, которую получим из матрицы $\mathbf{A}$ заменой ее $i$-й строки строкой $\mathbf{b}$

Рассмотренные обозначения могут быть скомбинированы. Например, запись

$$\mathbf{A}_{k.}^{ii}(\mathbf{b})$$

означает замену строки с номером $k$ матрицей $\mathbf{b}$ с последующим удалением строки с номером $i$ и столбца с номером $i$.

Как отмечалось в разделе [17]-2.2 мы можем определить два вида произведения матриц: либо произведение строк первой матрицы на столбцы второй, либо произведение столбцов первой матрицы на строки второй. Однако согласно теореме [17]-2.2.5 это произведение симметрично относительно операции транспонирования. Поэтому в статье мы ограничимся традиционным произведением строк первой матрицы на столбцы второй, и мы не будем указывать явно операцию как это было сделано в [17].

## 3. Предварительные замечания

Теорию определителей матриц с некоммутативными элементами, можно условно разделить на три группы относительно методов их определения. Обозначим через $\mathrm{M}(n,\mathbf{K})$ кольцо матриц с элементами из кольца $\mathbf{K}$. Один из способов определения детерминанта матрицы из $\mathrm{M}\,(n,\mathbf{K})$ — следующий ([1, 4, 5]).

**Определение 3.1.** Пусть функционал

$$\mathrm{d}:\mathrm{M}\,(n,\mathbf{K})\to\mathbf{K}$$

удовлетворяет следующим аксиомам.

**Аксиома 1.** $\mathrm{d}\,(\mathbf{A})=0$ тогда и только тогда, когда $\mathbf{A}$ – вырожденная (необратимая).
**Аксиома 2.** $\forall\mathbf{A},\mathbf{B}\in\mathrm{M}\,(n,\mathbf{K})$, $\mathrm{d}\,(\mathbf{A}\cdot\mathbf{B})=\mathrm{d}\,(\mathbf{A})\cdot\mathrm{d}\,(\mathbf{B})$.
**Аксиома 3.** Если матрица $\mathbf{A}'$ получается из матрицы $\mathbf{A}$ прибавлением ее произвольной строки умноженной слева с ее другой строкой или ее произвольного столбца умноженного справа с ее другим столбцом, тогда

$$\mathrm{d}\,(\mathbf{A}')=\mathrm{d}\,(\mathbf{A})$$

Тогда значение функционала d называется определителем матрицы $\mathbf{A}\in\mathrm{M}\,(n,\mathbf{K})$. □

Примерами таких функционалов являются известные определители Дьедонне и Стади. В [1] доказано, что определители, которые удовлетворяют аксиомы 1, 2, 3, принимают свое значение в некотором коммутативном подмножестве кольца. Также для них не имеет смысла такое свойство обычных определителей как разложение по минорам строки или столбца. Поэтому детерминантное представление обратной матрицы с помощью только таких определителей невозможно. Всё это является причиной, которая заставляет вводить детерминантные функционалы, не удовлетворяющие всем вышеприведенным аксиомам. При этом аксиома 1 рассматривается в [5] как необходимая для определения детерминанта.

При другом подходе определитель квадратной матрицы над некоммутативным кольцом строится как рациональная функция от элементов матрицы. Наибольшего успеха здесь достигли И. М. Гельфанд и В. С. Ретах своей теорией квазидетерминантов ([9, 10]). Введение в теорию квазидетерминантов будет изложено в разделе 6.

При третьем подходе детерминант матрицы с некоммутативными элементами определяется, как альтернированная сумма $n!$ произведений элементов матрицы, но с определенным фиксированным порядком множителей в каждом из них. Е. Г. Мур был первый, кто достиг выполнения ключевой аксиомы 1 при таком определении некоммутативного детерминанта. Мур сделал это не для всех квадратных матриц, а только для эрмитовых. Он определил детерминант эрмитовой матрицы[1] $\mathbf{A} = (a_{ij})_{n\times n}$ над телом с инволюцией индукцией по $n$ следующим образом (см. [5])

$$\mathrm{Mdet}\mathbf{A} = \begin{cases} a_{11}, & n = 1 \\ \sum\limits_{j=1}^{n} \varepsilon_{ij} a_{ij} \mathrm{Mdet}\left(\mathbf{A}(i \to j)\right), & n > 1 \end{cases} \tag{3.1}$$

Здесь $\varepsilon_{kj} = \begin{cases} 1, & i = j \\ -1, & i \neq j \end{cases}$, а через $\mathbf{A}(i \to j)$ обозначена матрица, которая получается из матрицы $\mathbf{A}$ последовательным применением замены $j$-го столбца $i$-м и вычеркивания $i$-ых строки и столбца. Другое определение этого детерминанта представлено в [1] в терминах подстановок:

$$\mathrm{Mdet}\,\mathbf{A} = \sum_{\sigma\in S_n} |\sigma| a_{n_{11}n_{12}} \cdot \ldots \cdot a_{n_{1l_1}n_{11}} \cdot a_{n_{21}n_{22}} \cdot \ldots \cdot a_{n_{rl_1}n_{r1}}.$$

Здесь $S_n$ – симметрическая группа $n$ элементов. Разложение подстановки $\sigma$ в циклы имеет вид:

$$\sigma = (n_{11} \ldots n_{1l_1})\,(n_{21} \ldots n_{2l_2}) \ldots (n_{r1} \ldots n_{rl_r})\,.$$

Однако, не существовало никакого обобщения определения детерминанта Мура для произвольных квадратных матриц. Ф. Д. Дайсон отмечал важность этой задачи в работе [5].

[1]**Эрмитова матрица** - это такая матрица $\mathbf{A} = (a_{ij})$, что $a_{ij} = \overline{a_{ji}}$.

В статьях [2, 3] Чен Лонгхуан предложил следующее определение детерминанта квадратной матрицы над телом кватернионов $\mathbf{A} = (a_{ij}) \in \mathrm{M}(n, \mathbf{H})$:

$$\det \mathbf{A} = \sum_{\sigma \in S_n} \varepsilon(\sigma)\, a_{n_1 i_2} \cdot a_{i_2 i_3} \cdot \ldots \cdot a_{i_s n_1} \cdot \ldots \cdot a_{n_r k_2} \cdot \ldots \cdot a_{k_l n_r},$$

$$\sigma = (n_1 i_2 \ldots i_s) \ldots (n_r k_2 \ldots k_r),$$

$$n_1 > i_2, i_3, \ldots, i_s; \ldots, n_r > k_2, k_3, \ldots, k_l,$$

$$n = n_1 > n_2 > \ldots > n_r \geq 1.$$

Несмотря на то, что такой определитель не удовлетворяет аксиоме 1, Л. Чен получил детерминантное представление обратной матрицы. Но его определитель также нельзя разложить по минорам строки или столбца (за исключением $n$-й строки). Поэтому Л. Чен также не получил классическую присоединенную матрицу. Для $\mathbf{A} = (\alpha_1, \ldots, \alpha_m)$ над телом кватернионов $\mathbb{H}$, если $\|\mathbf{A}\| := \det(\mathbf{A}^*\mathbf{A}) \neq 0$, тогда $\exists \mathbf{A}^{-1} = (b_{jk})$, где

$$\overline{b_{jk}} = \frac{1}{\|\mathbf{A}\|}\omega_{kj}, \quad (j, k = 1, 2, \ldots, n),$$

$$\omega_{kj} = \det(\alpha_1 \ldots \alpha_{j-1}\alpha_n\alpha_{j+1} \ldots \alpha_{n-1}\delta_k)^* (\alpha_1 \ldots \alpha_{j-1}\alpha_n\alpha_{j+1} \ldots \alpha_{n-1}\alpha_j).$$

Здесь $\alpha_i$ – $i$-й столбец $\mathbf{A}$, $\delta_k$ – $n$-мерный столбец с 1 в $k$-й строке и 0 в других. Л. Чен определил $\|\mathbf{A}\| := \det(\mathbf{A}^*\mathbf{A})$, как двойной детерминант. Решение правой системы линейных уравнений

$$\sum\nolimits_{j=1}^{n} \alpha_j x_j = \beta$$

над $\mathbb{H}$, если $\|\mathbf{A}\| \neq 0$, представлено следующей формулой, которую автор называет крамеровской

$$x_j = \|\mathbf{A}\|^{-1}\overline{\mathbf{D}_j}, \ \ \forall j = \overline{1, n}$$

где

$$\mathbf{D}_j = \det \begin{pmatrix} \alpha_1^* \\ \vdots \\ \alpha_{j-1}^* \\ \alpha_n^* \\ \alpha_{j+1}^* \\ \vdots \\ \alpha_{n-1}^* \\ \beta^* \end{pmatrix} \begin{pmatrix} \alpha_1 & \ldots & \alpha_{j-1} & \alpha_n & \alpha_{j+1} & \ldots & \alpha_{n-1} & \alpha_j \end{pmatrix}.$$

Здесь $\alpha_i$ – $i$-й столбец матрицы $\mathbf{A}$, $\alpha_i^*$ — $i$-я строка матрицы $\mathbf{A}^*$, $\beta^*$ — $n$-мерная вектор-строка сопряженная с $\beta$.

В данной работе мы рассматриваем теорию строчных и столбцовых определителей, которая развивает классический подход к определению детерминанта матрицы, но с заранее определенным порядком множителей в каждом из слагаемых детерминанта.

## 4. Кватернионная алгебра

**Кватернионная алгебра** $\mathbb{H}(a,b)$ (она также обозначается как $\left(\frac{a,b}{\mathbb{F}}\right)$) - это четырехмерное векторное пространство над полем $\mathbb{F}$ с базисом $\{1,i,j,k\}$ и следующими правилами умножения:

$$\begin{aligned} i^2 &= a, \\ j^2 &= b, \\ ij &= k, \\ ji &= -k. \end{aligned}$$

Поле $\mathbb{F}$ является центром кватернионной алгебры $\mathbb{H}(a,b)$.

В алгебре $\mathbb{H}(a,b)$ определены следующие отображения.

- Квадратичная форма
$$\mathrm{n} : x \in \mathbb{H} \to \mathrm{n}(x) \in \mathbb{F}$$
такая что
$$\mathrm{n}(x \cdot y) = \mathrm{n}(x)\mathrm{n}(y) \quad x, y \in \mathbb{H}$$
называется **нормой** в кватернионной алгебре $\mathbb{H}$.
- Линейное отображение
$$\mathrm{t} : x = x^0 + x^1 i + x^2 j + x^3 k \in \mathbb{H} \to \mathrm{t}(x) = 2x^0 \in \mathbb{F}$$
называется **следом кватерниона**. След удовлетворяет **условию перестановочности следа**
$$\mathrm{t}\,(q \cdot p) = \mathrm{t}\,(p \cdot q)$$
Из теоремы [17]-10.3.3 следует
$$\mathrm{t}(x) = \frac{1}{2}(x - ixi - jxj - kxk) \tag{4.1}$$
- Линейное отображение
$$x \to \overline{x} = \mathrm{t}(x) - x \tag{4.2}$$
является **инволюцией**. Инволюция имеет следующие свойства
$$\begin{aligned} \overline{\overline{x}} &= x \\ \overline{x+y} &= \overline{x} + \overline{y} \\ \overline{x \cdot y} &= \overline{y} \cdot \overline{x} \end{aligned}$$
При этом кватернион $\overline{x}$ будем называть **сопряженным к кватерниону** $x \in \mathbb{H}$. Норма и инволюция удовлетворяют следующему условию
$$\mathrm{n}\,(\overline{q}) = \mathrm{n}(q)$$
След и инволюция удовлетворяют следующему условию
$$\mathrm{t}\,(\overline{x}) = \mathrm{t}(x)$$
Из равенств (4.1), (4.2) следует
$$\overline{x} = -\frac{1}{2}(x + ixi + jxj + kxk)$$

В зависимости от выбора поля $\mathbb{F}$ и элементов $a, b$ на множестве всех кватернионных алгебр возможны два случая [19]:

1. $\left(\frac{a,b}{\mathbb{F}}\right)$ является алгеброй с делением.

2. $\left(\frac{a,b}{\mathbb{F}}\right)$ изоморфна алгебре всех $2\times 2$ матриц над полем $\mathbb{F}$. В этом случае кватернионная алгебра - расщепляемая.

Рассмотрим некоторые не изоморфные кватернионные алгебры с делением.

1. $\left(\frac{a,b}{\mathbb{R}}\right)$, где $\mathbb{R}$ - поле действительных чисел, изоморфна классическому телу кватернионов $\mathbf{H}$, когда $a<0$ и $b<0$. В противном случае, $\left(\frac{a,b}{\mathbb{R}}\right)$ - расщепляемая.

2. Над полем рациональных чисел $\mathbb{Q}$ существует бесконечно много не изоморфных кватернионных алгебр с делением $\left(\frac{a,b}{\mathbb{Q}}\right)$ в зависимости от выбора $a<0$ и $b<0$.

3. Пусть $\mathbb{Q}_p$ - поле $p$-адичних чисел, где $p$ - простое число. Для каждого простого числа $p$ существует единственная кватернионная алгебра с делением.

## 5. Введение в теорию строчных и столбцовых определителей

Строчные и столбцовые определители вводятся для матриц над кватернионной алгеброй $\mathbb{H}$. Следующие определения в теории подстановок необходимы нам для введения строчно-столбцовых определителей.

**Определение 5.1.** Пусть $S_n$ – симметрическая группа на множестве $I_n = \{1, 2, ..., n\}$ ([6], с. 70, упражнения 9, 10). Будем говорить, что подстановка $\sigma \in S_n$ образует прямое произведение независимых циклов, если ее запись в обычной двухрядной форме соответствует ее разложению в независимые циклы

$$\sigma = \begin{pmatrix} n_{11} & n_{12} & \dots & n_{1l_1} & \dots & n_{r1} & n_{r2} & \dots & n_{rl_r} \\ n_{12} & n_{13} & \dots & n_{11} & \dots & n_{r2} & n_{r3} & \dots & n_{r1} \end{pmatrix} \tag{5.1}$$

□

**Определение 5.2.** Будем говорить, что представление (5.1) подстановки $\sigma$ произведением независимых циклов является **упорядоченным слева**, если элементы, которые замыкают каждый из ее независимых циклов, записываются первыми слева в каждом из циклов

$$\sigma = (n_{11}n_{12}\dots n_{1l_1})\ (n_{21}n_{22}\dots n_{2l_2})\dots(n_{r1}n_{r2}\dots n_{rl_r})\,.$$

□

**Определение 5.3.** Будем говорить, что представление (5.1) подстановки $\sigma$ произведением независимых циклов является **упорядоченным справа**, если элементы, которые замыкают каждый из ее независимых циклов, записываются первыми справа в каждом из циклов

$$\sigma = (n_{12}\dots n_{1l_1}n_{11})\ (n_{22}\dots n_{2l_2}n_{21})\dots(n_{r2}\dots n_{rl_r}n_{r1})\,.$$

□

**Определение 5.4.** Пусть $S_n$ – симметрическая группа на множестве $I_n = \{1, 2, ..., n\}$. Пусть $i \in \overline{1,n}$. **Строчным определителем** по $i$-й строке матрицы $\mathbf{A} = (a_{ij}) \in \mathrm{M}\,(n, \mathbb{H})$ будем называть выражение

$$\mathrm{rdet}_i \mathbf{A} = \sum_{\sigma \in S_n} (-1)^{n-r} a_{i\, i_{k_1}} a_{i_{k_1} i_{k_1+1}} \dots a_{i_{k_1+l_1} i} \dots a_{i_{k_r} i_{k_r+1}} \dots a_{i_{k_r+l_r} i_{k_r}}$$

Это выражение является альтернированной суммой $n!$ мономов, составленных из элементов матрицы $\mathbf{A}$, в каждом из которых подстановка $\sigma \in S_n$ индексов образует прямое произведение независимых циклов

$$\sigma = (i\, i_{k_1} i_{k_1+1} \dots i_{k_1+l_1})\, (i_{k_2} i_{k_2+1} \dots i_{k_2+l_2}) \dots (i_{k_r} i_{k_r+1} \dots i_{k_r+l_r})$$

В упорядоченном слева разложении подстановки $\sigma$ первый слева цикл начинается слева индексом $i$, а все следующие циклы удовлетворяют условиям:

$$i_{k_2} < i_{k_3} < \dots < i_{k_r}, \quad i_{k_t} < i_{k_t+s}, \quad \left(\forall t = \overline{2,r}\right), \quad \left(\forall s = \overline{1,l_t}\right).$$

□

**Определение 5.5.** Пусть $S_n$ – симметрическая группа на множестве $J_n = \{1, 2, ..., n\}$. Пусть $j \in \overline{1,n}$. **Столбцовым определителем** по $j$-му столбцу матрицы $\mathbf{A} = (a_{ij}) \in \mathrm{M}\,(n, \mathbb{H})$ будем называть выражение

$$\mathrm{cdet}_{\mathrm{j}}\, \mathbf{A} = \sum_{\tau \in S_n} (-1)^{n-r}\, a_{j_{k_r} j_{k_r+l_r}} \dots a_{j_{k_r+1} i_{k_r}} \dots a_{j\, j_{k_1+l_1}} \dots a_{j_{k_1+1} j_{k_1}} a_{j_{k_1} j}.$$

Это выражение является альтернированной суммой $n!$ мономов, составленных из элементов матрицы $\mathbf{A}$, в каждом из которых подстановка $\tau \in S_n$ индексов образует прямое произведение независимых циклов

$$\tau = (j_{k_r+l_r} \dots j_{k_r+1} j_{k_r}) \dots (j_{k_2+l_2} \dots j_{k_2+1} j_{k_2})\ (j_{k_1+l_1} \dots j_{k_1+1} j_{k_1} j)$$

В упорядоченном справа разложении подстановки $\tau$ первый справа цикл начинается справа индексом $j$, а все следующие циклы удовлетворяют условиям:

$$j_{k_2} < j_{k_3} < \dots < j_{k_r}, \quad j_{k_t} < j_{k_t+s}, \quad \left(\forall t = \overline{2,r}\right), \quad \left(\forall s = \overline{1,l_t}\right).$$

□

*Замечание* 5.6. Особенностью столбцовых определителей является то, что при непосредственном их вычислении множители каждого из мономов записываются справа налево. □

В леммах 5.7 и 5.8 рассмотрено рекуррентное определение столбцовых и строчных определителей. Это определение является аналогом разложения определителя по строкам и столбцам в коммутативном случае.

*Лемма* 5.7. Пусть $R_{i\,j}$ – **правое алгебраическое дополнение** элемента $a_{i\,j}$ матрицы $\mathbf{A} = (a_{ij}) \in \mathrm{M}\,(n, \mathbb{H})$, а именно

$$\mathrm{rdet}_{\mathrm{i}}\, \mathbf{A} = \sum_{\mathrm{j}=1}^{\mathrm{n}} \mathrm{a}_{\mathrm{i\,j}} \cdot \mathrm{R}_{\mathrm{i\,j}} \quad i = \overline{1,n}$$

Тогда

$$R_{i\,j} = \begin{cases} -\mathrm{rdet}_j\ (\mathbf{A}^{ii}_{.j}(\mathbf{a}_{.\,i})) & i \neq j \\ \ \ \mathrm{rdet}_k\ \mathbf{A}^{ii} & i = j \end{cases}$$

$$k = \begin{cases} 2 & i = 1 \\ 1 & i > 1 \end{cases}$$

где матрица $(\mathbf{A}^{ii}_{.j}(\mathbf{a}_{.\,i}))$ получается из $\mathbf{A}$ последовательным применением замены $j$-го столбца $i$-м и вычеркивания $i$-х строки и столбца. □

*Лемма* 5.8. Пусть $L_{i\,j}$ – **левое алгебраическое дополнение** элемента $a_{i\,j}$ матрицы $\mathbf{A} = (a_{ij}) \in \mathrm{M}\,(n, \mathbb{H})$, а именно

$$\mathrm{cdet_j}\,\mathbf{A} = \sum_{\mathrm{i=1}}^{\mathrm{n}} \mathrm{L_{i\,j}} \cdot \mathrm{a_{i\,j}} \quad j = \overline{1,n}$$

Тогда

$$L_{i\,j} = \begin{cases} -\mathrm{cdet}_i\ (\mathbf{A}^{jj}_{i.}(\mathbf{a}_{j\,.})) & i \neq j \\ \ \ \mathrm{cdet}_k\ \mathbf{A}^{jj} & i = j \end{cases}$$

$$k = \begin{cases} 2 & j = 1 \\ 1 & j > 1 \end{cases}$$

где матрица $(\mathbf{A}^{jj}_{i.}(\mathbf{a}_{j\,.}))$ получается из $\mathbf{A}$ последовательным применением замены $i$-й строки $j$-й и вычеркивания $j$-х строки и столбца. □

*Замечание* 5.9. Очевидно, что каждому моному любого определенного выше детерминанта квадратной матрицы отвечает моном любого другого детерминанта, строчного или столбцового, такой, что оба они имеют одинаковый знак, состоят из одних и тех же множителей, - элементов матрицы, и отличаются только порядком их размещения. Кроме того, если элементы матрицы коммутируют, то $\mathrm{rdet}_1\,\mathbf{A} = \ldots = \mathrm{rdet}_n\mathbf{A} = \mathrm{cdet}_1\,\mathbf{A} = \ldots = \mathrm{cdet}_n\mathbf{A}$. □

Рассмотрим основные свойства строчного и столбцового детерминанта произвольной квадратной матрицы над телом $\mathbb{H}$, доказательства которых непосредственно следуют из определений.

**Теорема 5.10.** *Если одна из строк (столбцов) матрицы* $\mathbf{A} \in \mathrm{M}\,(n, \mathbb{H})$ *состоит из нолей, тогда для всех* $i = \overline{1,n}$

$$\mathrm{rdet}_i\,\mathbf{A} = 0, \quad \mathrm{cdet}_i\,\mathbf{A} = 0.$$

**Теорема 5.11.** *Если* $i$*-я строка матрицы* $\mathbf{A} = (a_{ij}) \in \mathrm{M}\,(n, \mathbb{H})$ *умножается слева на любое* $b \in \mathbb{H}$*, тогда для всех* $i = \overline{1,n}$

$$\mathrm{rdet}_i\,\mathbf{A}_{i\,.}(b \cdot \mathbf{a}_{.\,i}) = b \cdot \mathrm{rdet}_i\,\mathbf{A}.$$

**Теорема 5.12.** *Если* $j$*-й столбец матрицы* $\mathbf{A} = (a_{ij}) \in \mathrm{M}\,(n, \mathbb{H})$ *умножается справа на любое* $b \in \mathbb{H}$*, тогда для всех* $j = \overline{1,n}$

$$\mathrm{cdet}_j\,\mathbf{A}_{.\,j}(\mathbf{a}_{j\,.} \cdot b) = \mathrm{cdet}_j\,\mathbf{A} \cdot b.$$

**Теорема 5.13.** *Если для* $\mathbf{A} = (a_{ij}) \in \mathrm{M}(n, \mathbb{H})$ *найдется такое* $k \in I_n$, *что* $a_{kj} = b_j + c_j$ *для всех* $j = \overline{1, n}$, *тогда для любого* $i = \overline{1, n}$

$$\begin{aligned} \mathrm{rdet}_i\, \mathbf{A} &= \mathrm{rdet}_i\, \mathbf{A}_{k\,.}\,(\mathbf{b}) + \mathrm{rdet}_i\, \mathbf{A}_{k\,.}\,(\mathbf{c})\,, \\ \mathrm{cdet}_i\, \mathbf{A} &= \mathrm{cdet}_i\, \mathbf{A}_{k\,.}\,(\mathbf{b}) + \mathrm{cdet}_i\, \mathbf{A}_{k\,.}\,(\mathbf{c})\,. \end{aligned}$$

*где* $\mathbf{b} = (b_1, \ldots, b_n)$, $\mathbf{c} = (c_1, \ldots, c_n)$.

**Теорема 5.14.** *Если для* $\mathbf{A} = (a_{ij}) \in \mathrm{M}(n, \mathbb{H})$ *найдется такое* $k \in I_n$, *что* $a_{i\,k} = b_i + c_i$ *для всех* $i = \overline{1, n}$, *тогда для любого* $j = \overline{1, n}$

$$\begin{aligned} \mathrm{rdet}_j\, \mathbf{A} &= \mathrm{rdet}_j\, \mathbf{A}_{.\,k}\,(\mathbf{b}) + \mathrm{rdet}_j\, \mathbf{A}_{.\,k}\,(\mathbf{c})\,, \\ \mathrm{cdet}_j\, \mathbf{A} &= \mathrm{cdet}_j\, \mathbf{A}_{.\,k}\,(\mathbf{b}) + \mathrm{cdet}_j \mathbf{A}_{.\,k}\,(\mathbf{c})\,. \end{aligned}$$

*где* $\mathbf{b} = (b_1, \ldots, b_n)^T$, $\mathbf{c} = (c_1, \ldots, c_n)^T$.

**Теорема 5.15.** *Пусть* $\mathbf{A}^*$ *– матрица, эрмитово сопряженная к* $\mathbf{A} \in \mathrm{M}(n, \mathbb{H})$, *тогда* $\mathrm{rdet}_i\, \mathbf{A}^* = \overline{\mathrm{cdet}_i\, \mathbf{A}}$ *для любого* $i = \overline{1, n}$. □

Следующая теорема имеет ключевое значение в теории строчных и столбцовых определителей.

**Теорема 5.16.** *Если* $\mathbf{A} = (a_{ij}) \in \mathrm{M}(n, \mathbb{H})$ *– эрмитовая матрица, тогда*

$$\mathrm{rdet}_1 \mathbf{A} = \ldots = \mathrm{rdet}_n \mathbf{A} = \mathrm{cdet}_1 \mathbf{A} = \ldots = \mathrm{cdet}_n \mathbf{A} \in \mathbb{F}.$$

□

*Замечание* 5.17. Поскольку все столбцовые и строчные определители эрмитовой матрицы над телом $\mathbb{H}$ равны между собой, то мы можем однозначно ввести понятие определителя эрмитовой матрицы $\mathbf{A} \in \mathrm{M}(n, \mathbb{H})$:

$$\det \mathbf{A} := \mathrm{rdet}_i\, \mathbf{A} = \mathrm{cdet}_i\, \mathbf{A}, \ \left(\forall i = \overline{1, n}\right).$$

□

*Замечание* 5.18. Представляя детерминант эрмитовой матрицы как строчный определитель по произвольной $i$-й строке по лемме 5.7 имеем

$$\begin{aligned} \det \mathbf{A} &= -\sum_{\sigma \in S_n} a_{i\,j} \cdot \mathrm{rdet}_j\,(\mathbf{A}^{ii}_{.j}(\mathbf{a}_{.\,i})) + a_{i\,i} \cdot \mathrm{rdet}_k\, \mathbf{A}^{ii} \\ k &= \begin{cases} 2 & i = 1 \\ 1 & i > 1 \end{cases} \end{aligned} \tag{5.2}$$

Сравнивая выражения (3.1) и (5.2) для эрмитовой матрицы $\mathbf{A} \in \mathrm{M}(n, \mathbf{H})$, с очевидностью получаем, что введенный строчный определитель для эрмитовой матрицы совпадает с детерминантом Мура, а строчные и столбцовые определители для произвольных матриц являются его обобщением на множестве произвольных квадратных матриц. □

Свойства определителя эрмитовой матрицы полностью рассматриваются в [14] посредством ее строчных и столбцовых определителей. Среди всех рассмотрим следующие.

**Теорема 5.19.** *Если* $i$*-ю строку эрмитовой матрицы* $\mathbf{A} \in \mathrm{M}(n, \mathbb{H})$ *заменить левой линейной комбинацией других ее строк*

$$\mathbf{a}_{i\,.} = c_1 \mathbf{a}_{i_1\,.} + \ldots + c_k \mathbf{a}_{i_k\,.}$$

*где* $c_l \in \mathbb{H}$ *для всех* $l = \overline{1,k}$ *и* $\{i, i_l\} \subset I_n$*, тогда для всех* $i = \overline{1,n}$

$$\mathrm{cdet}_i \mathbf{A}_{i\,.}\left(c_1 \cdot \mathbf{a}_{i_1\,.} + \ldots + c_k \cdot \mathbf{a}_{i_k\,.}\right) = \mathrm{rdet}_i \mathbf{A}_{i\,.}\left(c_1 \cdot \mathbf{a}_{i_1\,.} + \ldots + c_k \cdot \mathbf{a}_{i_k\,.}\right) = 0.$$

□

**Теорема 5.20.** *Если* $j$*-й столбец эрмитовой матрицы* $\mathbf{A} \in \mathrm{M}(n, \mathbb{H})$ *заменить правой линейной комбинацией других ее столбцов*

$$\mathbf{a}_{.\,j} = \mathbf{a}_{.\,j_1} c_1 + \ldots + \mathbf{a}_{.\,j_k} c_k$$

*где* $c_l \in \mathbb{H}$ *для всех* $l = \overline{1,k}$ *и* $\{j, j_l\} \subset J_n$ *, тогда для всех* $j = \overline{1,n}$

$$\mathrm{cdet}_j \mathbf{A}_{.\,j}\left(\mathbf{a}_{.\,j_1} \cdot c_1 + \ldots + \mathbf{a}_{.\,j_k} \cdot c_k\right) = \mathrm{rdet}_j \mathbf{A}_{.\,j}\left(\mathbf{a}_{.\,j_1} \cdot c_1 + \ldots + \mathbf{a}_{.\,j_k} \cdot c_k\right) = 0.$$

□

Следующая теорема о детерминантном представлении матрицы обратной к эрмитовой непосредственно следует из этих свойств.

**Теорема 5.21.** *Для эрмитовой невырожденной матрицы* $\mathbf{A} \in \mathrm{M}(n, \mathbb{H})$*, (* $\det \mathbf{A} \neq 0$*), существует единственная правая* $(R\mathbf{A})^{-1}$ *и единственная левая обратная матрица* $(L\mathbf{A})^{-1}$*, при этом* $(R\mathbf{A})^{-1} = (L\mathbf{A})^{-1} =: \mathbf{A}^{-1}$*, которые обладают следующими детерминантными представлениями:*

$$(R\mathbf{A})^{-1} = \frac{1}{\det \mathbf{A}} \begin{pmatrix} R_{11} & R_{21} & \cdots & R_{n1} \\ R_{12} & R_{22} & \cdots & R_{n2} \\ \cdots & \cdots & \cdots & \cdots \\ R_{1n} & R_{2n} & \cdots & R_{nn} \end{pmatrix}$$

$$(L\mathbf{A})^{-1} = \frac{1}{\det \mathbf{A}} \begin{pmatrix} L_{11} & L_{21} & \cdots & L_{n1} \\ L_{12} & L_{22} & \cdots & L_{n2} \\ \cdots & \cdots & \cdots & \cdots \\ L_{1n} & L_{2n} & \cdots & L_{nn} \end{pmatrix}.$$

*Здесь* $R_{ij}$*,* $L_{ij}$ *являются соответственно левым и правым алгебраическими дополнениями для всех* $i, j = \overline{1,n}$*.* □

Для того, чтобы получить детерминантное представление для произвольной обратной матрицы над телом $\mathbb{H}$ рассматриваются ее правая $\mathbf{A}\mathbf{A}^*$ и левая $\mathbf{A}^*\mathbf{A}$ соответственные эрмитовые матрицы.

**Теорема 5.22** ([14])**.** *Если произвольный столбец матрицы* $\mathbf{A} \in \mathbb{H}^{m \times n}$ *является правой линейной комбинацией ее других столбцов, или произвольная строка матрицы* $\mathbf{A}^*$ *является левой линейной комбинацией ее других строк, тогда* $\det \mathbf{A}^*\mathbf{A} = 0$. □

Поскольку главная подматрица эрмитовой матрицы эрмитова, то главный базисный минор может быть определен по аналогии с коммутативным случаем, как ненулевой определитель ее главной подматрицы максимального порядка. Также вводим понятие **ранга эрмитовой матрицы по главным минорам**, как максимальный порядок отличного от нуля главного минора. Следующая теорема устанавливает соответствие между ним и рангом матрицы определенным как максимальное количество линейно независимых столбцов справа или линейно независимых строк слева, которые и образуют базис.

**Теорема 5.23** ([14])**.** *Ранг по главным минорам эрмитовой матрицы* $\mathbf{A}^*\mathbf{A}$ *равняется ее рангу, а также рангу матрицы* $\mathbf{A} \in \mathbb{H}^{m\times n}$. □

**Теорема 5.24** ([14])**.** *Если* $\mathbf{A} \in \mathbb{H}^{m\times n}$, *тогда произвольный столбец матрицы* $\mathbf{A}$ *является правой линейной комбинацией ее базисных столбцов или произвольная строка матрицы* $\mathbf{A}$ *является левой линейной комбинацией ее базисных строк.* □

Отсюда следует критерий вырожденности соответственной эрмитовой матрицы.

**Теорема 5.25** ([14])**.** *Правая линейная независимость столбцов матрицы* $\mathbf{A} \in \mathbb{H}^{m\times n}$ *или левая линейная независимость строк матрицы* $\mathbf{A}^*$ *является необходимым и достаточным условием того, что*

$$\det \mathbf{A}^*\mathbf{A} \neq 0$$

**Теорема 5.26** ([14])**.** *Если* $\mathbf{A} \in \mathrm{M}\,(n, \mathbb{H})$, *тогда* $\det \mathbf{A}\mathbf{A}^* = \det \mathbf{A}^*\mathbf{A}$. □

**Пример 5.27.** Рассмотрим матрицу

$$\mathbf{A} = \begin{pmatrix} a_{11} & a_{12} \\ a_{21} & a_{22} \end{pmatrix}$$

Тогда

$$\mathbf{A}^* = \begin{pmatrix} \overline{a_{11}} & \overline{a_{21}} \\ \overline{a_{12}} & \overline{a_{22}} \end{pmatrix}$$

Соответственно

$$\mathbf{A}\mathbf{A}^* = \begin{pmatrix} a_{11}\overline{a_{11}} + a_{12}\overline{a_{12}} & a_{11}\overline{a_{21}} + a_{12}\overline{a_{22}} \\ a_{21}\overline{a_{11}} + a_{22}\overline{a_{12}} & a_{21}\overline{a_{21}} + a_{22}\overline{a_{22}} \end{pmatrix}$$

$$\mathbf{A}^*\mathbf{A} = \begin{pmatrix} \overline{a_{11}}a_{11} + \overline{a_{21}}a_{21} & \overline{a_{11}}a_{12} + \overline{a_{21}}a_{22} \\ \overline{a_{12}}a_{11} + \overline{a_{22}}a_{21} & \overline{a_{12}}a_{12} + \overline{a_{22}}a_{22} \end{pmatrix}$$

Чтобы оценить детерминанты полученных эрмитовых матриц, рассмотрим например первый строчный определитель каждой матрицы. Согласно теореме 5.16 и замечанию 5.17

$$\det \mathbf{A}\mathbf{A}^* = \mathrm{rdet}_1 \mathbf{A}\mathbf{A}^*$$
$$\det \mathbf{A}^*\mathbf{A} = \mathrm{rdet}_1 \mathbf{A}^*\mathbf{A}$$

Согласно лемме 5.7

$$
\begin{aligned}
\det \mathbf{AA}^* &= (\mathbf{AA}^*)_{11}(\mathbf{AA}^*)_{22} - (\mathbf{AA}^*)_{12}(\mathbf{AA}^*)_{21}\\
&= (a_{11}\overline{a_{11}} + a_{12}\overline{a_{12}})(a_{21}\overline{a_{21}} + a_{22}\overline{a_{22}})\\
&-(a_{11}\overline{a_{21}} + a_{12}\overline{a_{22}})(a_{21}\overline{a_{11}} + a_{22}\overline{a_{12}})\\
&= a_{11}\overline{a_{11}}a_{21}\overline{a_{21}} + a_{12}\overline{a_{12}}a_{21}\overline{a_{21}}\\
&+a_{11}\overline{a_{11}}a_{22}\overline{a_{22}} + a_{12}\overline{a_{12}}a_{22}\overline{a_{22}}\\
&-a_{11}\overline{a_{21}}a_{21}\overline{a_{11}} - a_{12}\overline{a_{22}}a_{21}\overline{a_{11}}\\
&-a_{11}\overline{a_{21}}a_{22}\overline{a_{12}} - a_{12}\overline{a_{22}}a_{22}\overline{a_{12}}\\
&= a_{12}\overline{a_{12}}a_{21}\overline{a_{21}} + a_{11}\overline{a_{11}}a_{22}\overline{a_{22}}\\
&-a_{12}\overline{a_{22}}a_{21}\overline{a_{11}} - a_{11}\overline{a_{21}}a_{22}\overline{a_{12}}
\end{aligned}
\tag{5.3}
$$

$$
\begin{aligned}
\det \mathbf{A}^*\mathbf{A} &= (\mathbf{A}^*\mathbf{A})_{11}(\mathbf{A}^*\mathbf{A})_{22} - (\mathbf{A}^*\mathbf{A})_{12}(\mathbf{A}^*\mathbf{A})_{21}\\
&= (\overline{a_{11}}a_{11} + \overline{a_{21}}a_{21})(\overline{a_{12}}a_{12} + \overline{a_{22}}a_{22})\\
&-(\overline{a_{11}}a_{12} + \overline{a_{21}}a_{22})(\overline{a_{12}}a_{11} + \overline{a_{22}}a_{21})\\
&= \overline{a_{11}}a_{11}\overline{a_{12}}a_{12} + \overline{a_{21}}a_{21}\overline{a_{12}}a_{12}\\
&+\overline{a_{11}}a_{11}\overline{a_{22}}a_{22} + \overline{a_{21}}a_{21}\overline{a_{22}}a_{22}\\
&-\overline{a_{11}}a_{12}\overline{a_{12}}a_{11} - \overline{a_{21}}a_{22}\overline{a_{12}}a_{11}\\
&-\overline{a_{11}}a_{12}\overline{a_{22}}a_{21} - \overline{a_{21}}a_{22}\overline{a_{22}}a_{21}\\
&= \overline{a_{21}}a_{21}\overline{a_{12}}a_{12} + \overline{a_{11}}a_{11}\overline{a_{22}}a_{22}\\
&-\overline{a_{21}}a_{22}\overline{a_{12}}a_{11} - \overline{a_{11}}a_{12}\overline{a_{22}}a_{21}
\end{aligned}
\tag{5.4}
$$

Положительные слагаемые в равенствах (5.3), (5.4) являются действительными числами и, очевидно, совпадают. Чтобы доказать равенство

$$
a_{12}\overline{a_{22}}a_{21}\overline{a_{11}} + a_{11}\overline{a_{21}}a_{22}\overline{a_{12}} = \overline{a_{21}}a_{22}\overline{a_{12}}a_{11} + \overline{a_{11}}a_{12}\overline{a_{22}}a_{21}
\tag{5.5}
$$

используем свойство перестановочности следа элементов кватернионной алгебры, $\mathrm{t}(pq) = \mathrm{t}(qp)$. Действительно,

$$
a_{12}\overline{a_{22}}a_{21}\overline{a_{11}} + a_{11}\overline{a_{21}}a_{22}\overline{a_{12}} = a_{12}\overline{a_{22}}a_{21}\overline{a_{11}} + \overline{a_{12}\overline{a_{22}}a_{21}\overline{a_{11}}} = \mathrm{t}(a_{12}\overline{a_{22}}a_{21}\overline{a_{11}}),
$$

$$
\overline{a_{21}}a_{22}\overline{a_{12}}a_{11} + \overline{a_{11}}a_{12}\overline{a_{22}}a_{21} = \overline{\overline{a_{11}}a_{12}\overline{a_{22}}a_{21}} + \overline{a_{11}}a_{12}\overline{a_{22}}a_{21} = \mathrm{t}(\overline{a_{11}}a_{12}\overline{a_{22}}a_{21})
$$

Тогда из свойства перестановочности следа, следует (5.5). □

На основании теоремы 5.26 вводится понятие двойного определителя. Впервые это понятие было введено Л. Ченом ([2]).

**Определение 5.28.** Определитель эрмитовой матрицы $\mathbf{AA}^*$ называется **двойным определителем** матрицы $\mathbf{A} \in \mathrm{M}(n, \mathbb{H})$

$$
\mathrm{ddet}\mathbf{A} := \det(\mathbf{A}^*\mathbf{A}) = \det(\mathbf{AA}^*)
$$

□

В случае если $\mathbb{H}$ - классическое тело кватернионов $\mathbf{H}$ над полем действительных чисел, тогда имеет место теорема, которая устанавливает справедливость Аксиомы 1 для двойного определителя.

**Теорема 5.29.** *Если* $\{\mathbf{A}, \mathbf{B}\} \subset \mathrm{M}(n, \mathbf{H})$, *тогда* $\mathrm{ddet}(\mathbf{A} \cdot \mathbf{B}) = \mathrm{ddet}\mathbf{A} \cdot \mathrm{ddet}\mathbf{B}$. □

В общем случае не существует взаимосвязей между вырожденностью двойного определителя матрицы и её строчными и столбцовыми определителями. Продемонстрируем это в следующем примере.

**Пример 5.30.** Рассмотрим матрицу

$$A = \begin{pmatrix} i & j \\ j & -i \end{pmatrix}$$

Ее вторая строка получена из первой умножением слева на $k$. Тогда согласно теореме 5.25 $\mathrm{ddet} A = 0$. Действительно,

$$A^*A = \begin{pmatrix} -i & -j \\ -j & i \end{pmatrix} \cdot \begin{pmatrix} i & j \\ j & -i \end{pmatrix} = \begin{pmatrix} 2 & -2k \\ 2k & 2 \end{pmatrix}.$$

Тогда $\mathrm{ddet} A = 4 + 4k^2 = 0$, но

$$\mathrm{cdet}_1 A = \mathrm{cdet}_2 A = \mathrm{rdet}_1 A = \mathrm{rdet}_2 A = -i^2 - j^2 = 2$$

В то же время $rank A = 1$, что соответствует теореме 5.23. □

Установлено соответствие между двойным определителем и некоммутативными определителями Мура, Стади и Дьедонне,

$$\mathrm{ddet}\mathbf{A} = \mathrm{Mdet}\left(\mathbf{A}^*\mathbf{A}\right) = \mathrm{Sdet}\mathbf{A} = \mathrm{Ddet}^2\mathbf{A}$$

**Определение 5.31.** Пусть

$$\mathrm{ddet}\mathbf{A} = \mathrm{cdet}_j\left(\mathbf{A}^*\mathbf{A}\right) = \sum_i \mathbb{L}_{ij} \cdot a_{ij} \quad \forall j = \overline{1, n}$$

тогда будем называть $\mathbb{L}_{ij}$ **левым двойным алгебраическим дополнением** элемента $a_{ij}$ матрицы $\mathbf{A} \in \mathrm{M}\left(n, \mathbb{H}\right)$. □

**Определение 5.32.** Пусть

$$\mathrm{ddet}\mathbf{A} = \mathrm{rdet}_i\left(\mathbf{A}\mathbf{A}^*\right) = \sum_j a_{ij} \cdot \mathbb{R}_{ij} \quad \forall i = \overline{1, n}$$

тогда будем называть $\mathbb{R}_{ij}$ **правым двойным алгебраическим дополнением** элемента $a_{ij}$ матрицы $\mathbf{A} \in \mathrm{M}\left(n, \mathbf{H}\right)$. □

**Теорема 5.33.** *Необходимым и достаточным условием обратимости матрицы* $\mathbf{A} = (a_{ij}) \in \mathrm{M}(n, \mathbb{H})$ *является* $\mathrm{ddet}\mathbf{A} \neq 0$. *Тогда* $\exists \mathbf{A}^{-1} = \left(L\mathbf{A}\right)^{-1} = \left(R\mathbf{A}\right)^{-1}$, *где*

$$(5.6) \qquad \left(L\mathbf{A}\right)^{-1} = \left(\mathbf{A}^*\mathbf{A}\right)^{-1}\mathbf{A}^* = \frac{1}{\mathrm{ddet}\mathbf{A}} \begin{pmatrix} \mathbb{L}_{11} & \mathbb{L}_{21} & \dots & \mathbb{L}_{n1} \\ \mathbb{L}_{12} & \mathbb{L}_{22} & \dots & \mathbb{L}_{n2} \\ \dots & \dots & \dots & \dots \\ \mathbb{L}_{1n} & \mathbb{L}_{2n} & \dots & \mathbb{L}_{nn} \end{pmatrix}$$

$$(5.7) \qquad \left(R\mathbf{A}\right)^{-1} = \mathbf{A}^*\left(\mathbf{A}\mathbf{A}^*\right)^{-1} = \frac{1}{\mathrm{ddet}\mathbf{A}^*} \begin{pmatrix} \mathbb{R}_{11} & \mathbb{R}_{21} & \dots & \mathbb{R}_{n1} \\ \mathbb{R}_{12} & \mathbb{R}_{22} & \dots & \mathbb{R}_{n2} \\ \dots & \dots & \dots & \dots \\ \mathbb{R}_{1n} & \mathbb{R}_{2n} & \dots & \mathbb{R}_{nn} \end{pmatrix}$$

*и* $\mathbb{L}_{ij} = \mathrm{cdet}_j(\mathbf{A}^*\mathbf{A})_{.j}\left(\mathbf{A}^*_{.i}\right)$, $\mathbb{R}_{ij} = \mathrm{rdet}_i(\mathbf{A}\mathbf{A}^*)_{i.}\left(\mathbf{A}^*_{j.}\right)$, $\left(\forall i, j = \overline{1, n}\right)$. □

**Следствие 5.34.** *Пусть* $\mathrm{ddet}\mathbf{A} \neq 0$. *Тогда*

$$\mathbb{L}_{ij} = \mathbb{R}_{ij}$$

□

*Замечание* 5.35. Теорема 5.33 представляет обратную матрицу $\mathbf{A}^{-1}$ произвольной квадратной матрицы $\mathbf{A} \in \mathrm{M}(n, \mathbb{H})$, если $\mathrm{ddet}\mathbf{A} \neq 0$, через аналог классической присоединенной. Если мы обозначим ее $\mathrm{Adj}[[\mathbf{A}]]$, тогда над телом $\mathbb{H}$ имеет место следующая формула:

$$\mathbf{A}^{-1} = \frac{\mathrm{Adj}[[\mathbf{A}]]}{\mathrm{ddet}\mathbf{A}}.$$

□

Очевидным следствием детерминантного представления обратной матрицы через аналог классической присоединенной матрицы является правило Крамера.

**Теорема 5.36.** *Пусть*

$$\mathbf{A} \cdot \mathbf{x} = \mathbf{y} \tag{5.8}$$

*правая система линейных уравнений с матрицей коэффициентов* $\mathbf{A} \in \mathrm{M}(n, \mathbb{H})$, *столбцом свободных элементов* $\mathbf{y} = (y_1, \ldots, y_n)^T \in \mathbb{H}^{n \times 1}$ *и столбцом неизвестных* $\mathbf{x} = (x_1, \ldots, x_n)^T$. *Если* $\mathrm{ddet}\mathbf{A} \neq 0$, *тогда система линейных уравнений* (5.8) *имеет единственное решение, которое представляется формулой:*

$$x_j = \frac{\mathrm{cdet}_j (\mathbf{A}^* \mathbf{A})_{.j} (\mathbf{f})}{\mathrm{ddet}\mathbf{A}} \qquad \forall j = \overline{1, n} \tag{5.9}$$

*где* $\mathbf{f} = \mathbf{A}^* \mathbf{y}$. □

**Теорема 5.37.** *Пусть*

$$\mathbf{x} \cdot \mathbf{A} = \mathbf{y} \tag{5.10}$$

*– левая система линейных уравнений с матрицей коэффициентов* $\mathbf{A} \in \mathrm{M}(n, \mathbb{H})$, *строкой свободных элементов* $\mathbf{y} = (y_1, \ldots, y_n) \in \mathbb{H}^{1 \times n}$ *и строкой неизвестных* $\mathbf{x} = (x_1, \ldots, x_n)$. *Если* $\mathrm{ddet}\mathbf{A} \neq 0$, *тогда система линейных уравнений* (5.10) *имеет единственное решение, которое представляется формулой:*

$$x_i = \frac{\mathrm{rdet}_i \left(\mathbf{A}\mathbf{A}^*\right)_{i.} (\mathbf{z})}{\mathrm{ddet}\mathbf{A}} \qquad \forall i = \overline{1, n} \tag{5.11}$$

*где* $\mathbf{z} = \mathbf{y}\mathbf{A}^*$. □

Формулы (5.9) и (5.11) являются естественным и очевидным обобщением правила Крамера для квадратных систем линейных уравнений над кватернионной алгеброй с делением. Еще более близкую аналогию можно получить в следующих частных случаях, которые следуют из теоремы 5.21.

**Теорема 5.38.** *Пусть в правой системе линейных уравнений* (5.8) *матрица коэффициентов* $\mathbf{A} \in \mathrm{M}(n, \mathbb{H})$ *– эрмитова. Тогда система имеет единственное решение, которое представляется формулой:*

$$x_j = \frac{\mathrm{cdet}_j \mathbf{A}_{.j} (\mathbf{y})}{\det \mathbf{A}}, \quad \left(\forall j = \overline{1, n}\right).$$

□

**Теорема 5.39.** *Пусть в левой системе линейных уравнений (5.10) матрица коэффициентов* $\mathbf{A} \in \mathrm{M}(n, \mathbb{H})$ *– эрмитова. Тогда система имеет единственное решение, которое представляется формулой:*

$$x_i = \frac{\mathrm{rdet}_i \mathbf{A}_{i.}(\mathbf{y})}{\det \mathbf{A}}, \quad \left(\forall i = \overline{1, n}\right).$$

□

В рамках теории строчных-столбцовых определителей получено также правило Крамера для правого $\mathbf{AX} = \mathbf{B}$, левого $\mathbf{XA} = \mathbf{B}$ и двустороннего $\mathbf{AXB} = \mathbf{C}$ матричных уравнений ([15]). Также получено детерминантное представление Мура-Пенроуза обратной матрицы над телом кватернионов и правило Крамера для нормального решения правой и левой систем линейных уравнений ([16]).

## 6. Квазидетерминант

**Теорема 6.1.** *Предположим, что матрица*[2]

$$\mathbf{A} = \begin{pmatrix} a_{11} & ... & a_{1n} \\ ... & ... & ... \\ a_{n1} & ... & a_{nn} \end{pmatrix}$$

*имеет обратную матрицу* $A^{-1}$. *Тогда минор обратной матрицы удовлетворяет следующему равенству, при условии, что рассматриваемые обратные матрицы существуют,*

$$((\mathbf{A}^{-1})_{IJ})^{-1} = \mathbf{A}_{JI} - \mathbf{A}_{J.}^{.I}(\mathbf{A}^{JI})^{-1}\mathbf{A}_{.I}^{J.} \tag{6.1}$$

*Доказательство.* Определение обратной матрицы приводит к системе линейных уравнений

$$\mathbf{A}^{JI}(\mathbf{A}^{-1})_{.J}^{I.} + \mathbf{A}_{.I}^{J.}(\mathbf{A}^{-1})_{IJ} = 0 \tag{6.2}$$

$$\mathbf{A}_{J.}^{.I}(\mathbf{A}^{-1})_{.J}^{I.} + \mathbf{A}_{JI}(\mathbf{A}^{-1})_{IJ} = E \tag{6.3}$$

Мы умножим (6.2) на $\left(\mathbf{A}^{JI}\right)^{-1}$

$$(\mathbf{A}^{-1})_{.J}^{I.} + (\mathbf{A}^{JI})^{-1}\mathbf{A}_{.I}^{J.}(\mathbf{A}^{-1})_{IJ} = 0 \tag{6.4}$$

Теперь мы можем подставить (6.4) в (6.3)

$$\mathbf{A}_{JI}(\mathbf{A}^{-1})_{IJ} - \mathbf{A}_{J.}^{.I}(\mathbf{A}^{JI})^{-1}\mathbf{A}_{.I}^{J.}(\mathbf{A}^{-1})_{IJ} = E \tag{6.5}$$

(6.1) следует из (6.5). □

**Следствие 6.2.** *Предположим, что матрица* $\mathbf{A}$ *имеет обратную матрицу. Тогда элементы обратной матрицы удовлетворяют равенству*

$$((\mathbf{A}^{-1})_{ij})^{-1} = a_{ji} - \mathbf{A}_{j.}^{.i}(\mathbf{A}^{ji})^{-1}\mathbf{A}_{.i}^{j.} \tag{6.6}$$

□

[2]Это утверждение и его доказательство основаны на утверждении 1.2.1 из [7] для матриц над свободным кольцом с делением.

**Пример 6.3.** Рассмотрим матрицу

$$\mathbf{A} = \begin{pmatrix} a_{11} & a_{12} \\ a_{21} & a_{22} \end{pmatrix}$$

Согласно (6.6)

$$(A^{-1})_{11} = (a_{11} - a_{12}(a_{22})^{-1} \ a_{21})^{-1} \tag{6.7}$$

$$(A^{-1})_{21} = (a_{21} - a_{22}(a_{12})^{-1} \ a_{11})^{-1} \tag{6.8}$$

$$(A^{-1})_{12} = (a_{12} - a_{11}(a_{21})^{-1} \ a_{22})^{-1} \tag{6.9}$$

$$(A^{-1})_{22} = (a_{22} - a_{21}(a_{11})^{-1} \ a_{12})^{-1} \tag{6.10}$$

□

Мы называем матрицу

$$\mathcal{H}A = ((\mathcal{H}A)_{ij}) = ((a_{ji})^{-1}) \tag{6.11}$$

**обращением Адамара матрицы** [3]

$$\mathbf{A} = \begin{pmatrix} a_{11} & ... & a_{1n} \\ ... & ... & ... \\ a_{n1} & ... & a_{nn} \end{pmatrix}$$

**Определение 6.4.** $(ji)$**-квазидетерминант** матрицы $\mathbf{A}$ - это формальное выражение

$$|\mathbf{A}|_{ji} = (\mathcal{H}\mathbf{A}^{-1})_{ji} = ((\mathbf{A}^{-1})_{ij})^{-1} \tag{6.12}$$

Мы можем рассматривать $(ji)$-квазидетерминант как элемент матрицы $|A|$, которую мы будем называть **квазидетерминантом**. □

**Теорема 6.5.** *Выражение для* $(ji)$*-квазидетерминанта имеет форму*

$$|\mathbf{A}|_{ji} = a_{ji} - \mathbf{A}_{j\cdot}^{\cdot i}(\mathbf{A}^{ji})^{-1}\mathbf{A}_{\cdot i}^{j\cdot} \tag{6.13}$$

$$|\mathbf{A}|_{ji} = a_{ji} - \mathbf{A}_{j\cdot}^{\cdot i}\,\mathcal{H}|\mathbf{A}^{ji}|\,\mathbf{A}_{\cdot i}^{j\cdot} \tag{6.14}$$

*Доказательство.* Утверждение следует из (6.6) и (6.12). □

**Теорема 6.6.** *Пусть*

$$A = \begin{pmatrix} 1 & 0 \\ 0 & 1 \end{pmatrix} \tag{6.15}$$

*Тогда*

$$A^{-1} = \begin{pmatrix} 1 & 0 \\ 0 & 1 \end{pmatrix} \tag{6.16}$$

*Доказательство.* Из (6.7) и (6.10) очевидно, что $(A^{-1})_{11} = 1$ и $(A^{-1})_{22} = 1$. Тем не менее выражение для $(A^{-1})_{21}$ и $(A^{-1})_{12}$ не может быть определено из (6.8) и (6.9) так как $(a_{21}-a_{22}(a_{12})^{-1}\ a_{11})^{-1}=(a_{12}-a_{11}(a_{21})^{-1}\ a_{22})^{-1}=0$. Мы можем преобразовать эти выражения. Например

$$\begin{aligned}(A^{-1})_{21}&=(a_{21}-a_{22}(a_{12})^{-1}\ a_{11})^{-1}\\&=(a_{11}((a_{11})^{-1}\ a_{12}-(a_{21})^{-1}\ a_{22}))^{-1}\\&=((a_{21})^{-1}\ a_{11}(a_{21}(a_{11})^{-1}\ a_{12}-a_{22}))^{-1}\\&=(a_{11}(a_{21}(a_{11})^{-1}\ a_{12}-a_{22}))^{-1}\ a_{21}\end{aligned}$$

[3] [8]-стр. 4

Мы непосредственно видим, что $(A^{-1})_{21}=0$. Таким же образом мы можем найти, что $(A^{-1})_{12}=0$. Это завершает доказательство (6.16). □

Из доказательства теоремы 6.6 видно, что мы не всегда можем пользоваться равенством (6.6) для определения элементов обратной матрицы и необходимы дополнительные преобразования для решения этой задачи. Из теоремы [17]-4.6.3 следует, что если

$$\operatorname{rank}\begin{pmatrix}a_{11}&...&a_{1n}\\...&...&...\\a_{n1}&...&a_{nn}\end{pmatrix}\le n-2$$

то $|\mathbf{A}|_{ij}$, $i=1,\ ...,\ n$, $j=1,\ ...,\ n$, не определён. Из этого следует, что хотя квазидетерминант является мощным инструментом, возможность пользоваться определителем является серьёзным преимуществом.

**Теорема 6.7.** *Если матрица* **A** *имеет обратную матрицу, то для любых матриц* **B** *и* **C** *из равенства*

$$\mathbf{BA}=\mathbf{CA}\tag{6.17}$$

*следует равенство*

$$\mathbf{B}=\mathbf{C}\tag{6.18}$$

*Доказательство.* Равенство (6.18) следует из (6.17), если обе части равенства (6.17) умножить на $\mathbf{A}^{-1}$. □

**Теорема 6.8.** *Решение невырожденной системы линейных уравнений*

$$\mathbf{A}x=b\tag{6.19}$$

*определено однозначно и может быть записано в любой из следующих форм*[4]

$$x=\mathbf{A}^{-1}b\tag{6.20}$$

$$x=\mathcal{H}|\mathbf{A}|\,b\tag{6.21}$$

*Доказательство.* Умножая обе части равенства (6.19) слева на $\mathbf{A}^{-1}$, мы получим (6.20). Пользуясь определением 6.4, мы получим (6.21). Решение системы единственно в силу теоремы 6.7. □

[4]Смотри аналогичное утверждение в теореме [7]-1.6.1.

## 7. Соответствие между строчно-столбцовыми определителями и квазидетерминантами

**Теорема 7.1.** *Если* $\mathbf{A} \in \mathrm{M}(n, \mathbb{H})$ *- обратимая матрица, то для произвольных* $p, q = 1, \ldots, n$ *мы имеем следующие представления квазидетерминанта*

$$|\mathbf{A}|_{pq} = \mathrm{ddet}\mathbf{A}\ (\mathbb{L}_{\mathrm{pq}})^{-1} = \mathrm{ddet}\mathbf{A}\ (\mathrm{cdet_q}(\mathbf{A}^*\mathbf{A})_{.\,\mathrm{q}}(\mathbf{A}^*_{.\mathrm{p}}))^{-1} \tag{7.1}$$

$$= \frac{\mathrm{ddet}\mathbf{A}}{\mathrm{n}(\mathrm{cdet}_q(\mathbf{A}^*\mathbf{A})_{.\,q}(\mathbf{A}^*_{.p}))} \overline{\mathrm{cdet}_q(\mathbf{A}^*\mathbf{A})_{.\,q}(\mathbf{A}^*_{.p})} \tag{7.2}$$

$$|\mathbf{A}|_{pq} = \mathrm{ddet}\mathbf{A}\ (\mathbb{R}_{\mathrm{pq}})^{-1} = \mathrm{ddet}\mathbf{A}\ (\mathrm{rdet_p}(\mathbf{A}\mathbf{A}^*)_{\mathrm{p}\,.}(\mathbf{A}^*_{\mathrm{q}\,.}))^{-1} \tag{7.3}$$

$$= \frac{\mathrm{ddet}\mathbf{A}}{\mathrm{n}(\mathrm{rdet}_p(\mathbf{A}\mathbf{A}^*)_{p\,.}(\mathbf{A}^*_{q\,.}))} \overline{\mathrm{rdet}_p(\mathbf{A}\mathbf{A}^*)_{p\,.}(\mathbf{A}^*_{q\,.})} \tag{7.4}$$

*Доказательство.* Пусть $\mathbf{A}^{-1} = (b_{ij})$ - матрица обратная матрице $\mathbf{A}$. Равенство (6.12) раскрывает связь между квазидетерминантом $|\ \mathbf{A}\ |_{p,q}$ матрицы $\mathbf{A} \in \mathrm{M}(n, \mathbb{H})$ и элементами обратной матрицы $\mathbf{A}^{-1} = (b_{ij})$, а именно

$$|\ \mathbf{A}\ |_{pq} = b_{qp}^{-1}$$

для всех $p$, $q = 1, \ldots, n$. В тоже время теория строчных-столбцовых определителей (теорема 5.33) дает нам представления обратной матрицы через ее двойные левые (5.6) и правые (5.7) алгебраические дополнения. Таким образом, соответственно, получим

$$|\ \mathbf{A}\ |_{pq} = b_{qp}^{-1} = \left(\frac{\mathbb{L}_{pq}}{\mathrm{ddet}\mathbf{A}}\right)^{-1} = \left(\frac{\mathrm{cdet}_q(\mathbf{A}^*\mathbf{A})_{.\,q}\left(\mathbf{A}^*_{.p}\right)}{\mathrm{ddet}\mathbf{A}}\right)^{-1}, \tag{7.5}$$

$$|\ \mathbf{A}\ |_{pq} = b_{qp}^{-1} = \left(\frac{\mathbb{R}_{pq}}{\mathrm{ddet}\mathbf{A}}\right)^{-1} = \left(\frac{\mathrm{rdet}_p(\mathbf{A}\mathbf{A}^*)_{p\,.}\left(\mathbf{A}^*_{q\,.}\right)}{\mathrm{ddet}\mathbf{A}}\right)^{-1}. \tag{7.6}$$

Так как $\mathrm{ddet}\mathbf{A} \neq 0 \in \mathbb{F}$, то $\exists(\mathrm{ddet}\mathbf{A})^{-1} \in \mathbb{F}$. В свою очередь

$$\mathrm{cdet}_q(\mathbf{A}^*\mathbf{A})_{.\,q}\left(\mathbf{A}^*_{.p}\right)^{-1} = \frac{\overline{\mathrm{cdet}_q(\mathbf{A}^*\mathbf{A})_{.\,q}\left(\mathbf{A}^*_{.p}\right)}}{\mathrm{n}(\mathrm{cdet}_q(\mathbf{A}^*\mathbf{A})_{.\,q}\left(\mathbf{A}^*_{.p}\right))}, \tag{7.7}$$

$$\mathrm{rdet}_p(\mathbf{A}\mathbf{A}^*)_{p\,.}\left(\mathbf{A}^*_{q\,.}\right)^{-1} = \frac{\overline{\mathrm{rdet}_p(\mathbf{A}\mathbf{A}^*)_{p\,.}\left(\mathbf{A}^*_{q\,.}\right)}}{\mathrm{n}(\mathrm{rdet}_p(\mathbf{A}\mathbf{A}^*)_{p\,.}\left(\mathbf{A}^*_{q\,.}\right))}. \tag{7.8}$$

Подставив (7.7) в (7.5), а (7.8) в (7.6), соответственно получим (7.2) и (7.4).

Теорема доказана. □

Формула (7.2) дает явное представление квазидетерминанта $|\ \mathbf{A}\ |_{p,q}$ матрицы $\mathbf{A} \in \mathrm{M}(n, \mathbb{H})$ для всех $p, q = 1, \ldots, n$ через столбцовый определитель ее соответственной левой эрмитовой матрицы $\mathbf{A}^*\mathbf{A}$, а (7.4) - через строчный определитель ее соответственной правой эрмитовой $\mathbf{A}\mathbf{A}^*$.

**Пример 7.2.** Рассмотрим матрицу

$$\mathbf{A} = \begin{pmatrix} a_{11} & a_{12} \\ a_{21} & a_{22} \end{pmatrix}$$

Согласно (6.13)

$$|A| = \begin{pmatrix} a_{11} - a_{12}(a_{22})^{-1}\ a_{21} & a_{12} - a_{11}(a_{21})^{-1}\ a_{22} \\ a_{21} - a_{22}(a_{12})^{-1}\ a_{11} & a_{22} - a_{21}(a_{11})^{-1}\ a_{12} \end{pmatrix} \tag{7.9}$$

Наша задача найти этот квазидетерминант, пользуясь теоремой 7.1. Очевидно, что

$$\mathbf{A}^* = \begin{pmatrix} \overline{a_{11}} & \overline{a_{21}} \\ \overline{a_{12}} & \overline{a_{22}} \end{pmatrix} \quad \mathbf{A}^*\mathbf{A} = \begin{pmatrix} \mathrm{n}(a_{11}) + \mathrm{n}(a_{21}) & \overline{a_{11}}a_{12} + \overline{a_{21}}a_{22} \\ \overline{a_{12}}a_{11} + \overline{a_{22}}a_{21} & \mathrm{n}(a_{12}) + \mathrm{n}(a_{22}) \end{pmatrix}.$$

Вычисляем необходимые определители

$$\begin{aligned}
\mathrm{ddet}\mathbf{A} &= \mathrm{rdet}_1(\mathbf{A}^*\mathbf{A}) \\
&= (\mathrm{n}(a_{11}) + \mathrm{n}(a_{21})) \cdot (\mathrm{n}(a_{12}) + \mathrm{n}(a_{22})) \\
&\quad -(\overline{a_{11}}a_{12} + \overline{a_{21}}a_{22}) \cdot (\overline{a_{12}}a_{11} + \overline{a_{22}}a_{21}) \\
&= \mathrm{n}(a_{11})\mathrm{n}(a_{12}) + \mathrm{n}(a_{11})\mathrm{n}(a_{22}) + \mathrm{n}(a_{21})\mathrm{n}(a_{12}) + \mathrm{n}(a_{21})\mathrm{n}(a_{22}) \\
&\quad -\overline{a_{11}}a_{12}\overline{a_{12}}a_{11} - \overline{a_{11}}a_{12}\overline{a_{22}}a_{21} - \overline{a_{21}}a_{22}\overline{a_{12}}a_{11} - \overline{a_{21}}a_{22}\overline{a_{22}}a_{21} \\
&= \mathrm{n}(a_{11})\mathrm{n}(a_{22}) + \mathrm{n}(a_{21})\mathrm{n}(a_{12}) - (\overline{a_{11}}a_{12}\overline{a_{22}}a_{21} + \overline{\overline{a_{11}}a_{12}\overline{a_{22}}a_{21}}) \\
&= \mathrm{n}(a_{11})\mathrm{n}(a_{22}) + \mathrm{n}(a_{21})\mathrm{n}(a_{12}) - \mathrm{t}(\overline{a_{11}}a_{12}\overline{a_{22}}a_{21})
\end{aligned}$$

$$\begin{aligned}
\mathrm{cdet}_1(\mathbf{A}^*\mathbf{A})_{.1}(\mathbf{a}^*_{.2}) &= \mathrm{cdet}_1 \begin{pmatrix} \overline{a_{21}} & \overline{a_{11}}a_{12} + \overline{a_{21}}a_{22} \\ \overline{a_{22}} & \mathrm{n}(a_{12}) + \mathrm{n}(a_{22}) \end{pmatrix} \\
&= \mathrm{n}(a_{12})\overline{a_{21}} + \mathrm{n}(a_{22})\overline{a_{21}} - \overline{a_{11}}a_{12}\overline{a_{22}} - \overline{a_{21}}a_{22}\overline{a_{22}} \\
&= \mathrm{n}(a_{12})\overline{a_{21}} - \overline{a_{11}}a_{12}\overline{a_{22}}.
\end{aligned}$$

Тогда,

$$\begin{aligned}
\overline{\mathrm{cdet}_1(\mathbf{A}^*\mathbf{A})_{.1}(\mathbf{a}^*_{.2})} &= \mathrm{n}(a_{12})a_{21} - a_{22}\overline{a_{12}}a_{11}, \\
\mathrm{n}(\mathrm{cdet}_1(\mathbf{A}^*\mathbf{A})_{.1}(\mathbf{a}^*_{.2})) &= \overline{\mathrm{cdet}_1(\mathbf{A}^*\mathbf{A})_{.1}(\mathbf{a}^*_{.2})} \cdot \mathrm{cdet}_1(\mathbf{A}^*\mathbf{A})_{.1}(\mathbf{a}^*_{.2}) \\
&= (\mathrm{n}(a_{12})a_{21} - a_{22}\overline{a_{12}}a_{11}) \cdot (\mathrm{n}(a_{12})\overline{a_{21}} - \overline{a_{11}}a_{12}\overline{a_{22}}) \\
&= \mathrm{n}^2(a_{12})\mathrm{n}(a_{21}) - \mathrm{n}(a_{12})a_{21}\overline{a_{11}}a_{12}\overline{a_{22}} \\
&\quad -\mathrm{n}(a_{12})a_{22}\overline{a_{12}}a_{11}\overline{a_{21}} + a_{22}\overline{a_{12}}a_{11}\overline{a_{11}}a_{12}\overline{a_{22}} \\
&= \mathrm{n}(a_{12})(\mathrm{n}(a_{12})\mathrm{n}(a_{21}) - \mathrm{t}(\overline{a_{11}}a_{12}\overline{a_{22}}a_{21}) + \mathrm{n}(a_{21})\mathrm{n}(a_{12})) \\
&= \mathrm{n}(a_{12})\mathrm{ddet}\mathbf{A}.
\end{aligned}$$

Следуя формуле (7.2), получим

$$\begin{aligned}
|\mathbf{A}|_{21} &= \frac{\mathrm{ddet}\mathbf{A}}{\mathrm{n}(\mathrm{cdet}_1(\mathbf{A}^*\mathbf{A})_{.1}(\mathbf{a}^*_{.2}))}\overline{\mathrm{cdet}_1(\mathbf{A}^*\mathbf{A})_{.1}(\mathbf{a}^*_{.2})} \\
&= \frac{\mathrm{ddet}\mathbf{A}}{\mathrm{n}(a_{12})\mathrm{ddet}\mathbf{A}}\overline{\mathrm{cdet}_1(\mathbf{A}^*\mathbf{A})_{.1}(\mathbf{a}^*_{.2})} \\
&= \frac{1}{\mathrm{n}(a_{12})} \cdot \overline{\mathrm{cdet}_1(\mathbf{A}^*\mathbf{A})_{.1}(\mathbf{a}^*_{.2})} \\
&= \frac{1}{\mathrm{n}(a_{12})} \cdot (\mathrm{n}(a_{12})a_{21} - a_{22}\overline{a_{12}}a_{11}) \\
&= a_{21} - a_{22}(a_{12})^{-1}a_{11}.
\end{aligned} \tag{7.10}$$

Последнее выражение в (7.10) совпадает с выражением $|\mathbf{A}|_{21}$ в (7.9). □

## 8. Список литературы

# 9. Предметный указатель

## 10. Специальные символы и обозначения